\documentclass[a4paper,12pt]{article}
\usepackage[cp1251]{inputenc} 
\usepackage[english]{babel}
\usepackage{amsfonts,amssymb}
\title{The tensor functor from the category of $\A$-algebras into the category of differential modules
with $\infty$-simplicial faces.}
\author{S.V. Lapin}
\date{}
\newcommand{\F}{F_{\infty}}
\newcommand{\D}{D_{\infty}}
\newcommand{\bu}{\bullet}
\newcommand{\E}{E_{\infty}}
\newcommand{\A}{A_{\infty}}

\newcommand{\p}{\partial}

\begin{document}
\maketitle

\begin{abstract} In the present paper the tensor functor
from the category of $\A$-algebras into the category of
differential modules with $\infty$-simplicial faces is
constructed. Further, it is showed that this functor sends
homotopy equivalent $\A$-al\-geb\-ras into homotopy equivalent
differential modules with $\infty$-simplicial faces.
\end{abstract}

In \cite{Lap2}-\cite{Lap7} the homotopy technique of differential
modules with $\infty$-simplicial faces was developed. This
homotopy technique is closely related to the homotopy technique of
$\D$-differential modules that developed in
\cite{Lap9}-\cite{Lap17}. On the other hand, in \cite{Lap4} by any
$\A$-algebra the tensor differential module with
$\infty$-simplicial faces was constructed. By using the technique
of these tensor differential modules with $\infty$-simplicial
faces that defined by $\A$-algebras the concepts of cyclic
homology, dihedral homology and reflexive homology of
$\A$-algebras were developed in \cite{Lap18}-\cite{Lap20}. The
fact \cite{Lap4} that by each $\A$-algebra corresponds a tensor
differential module with $\infty$-simplicial faces give rise the
very important and interesting problem of extension of this
correspondence until a tensor functor from the category of
$\A$-algebras into the category of differential modules with
$\infty$-simplicial faces.

The present paper is devoted to solving the specified above
problem. The paper consists of three paragraphs. In the first
paragraph, we recall the necessary information from
\cite{Lap2}-\cite{Lap7} related to the notion of a differential
module with $\infty$-simplicial faces. In the second paragraph, we
recall the necessary information from \cite{Kad}-\cite{Smir}
related to the notion of an $\A$-algebra. Moreover, in this
paragraph, we pay special attention to the problem of correctly
writing signs in the structural relations of $\A$-algebras because
it is very important in the further considerations. In the third
paragraph, we construct a tensor functor from the category of
$\A$-algebras into the category of differential modules with
$\infty$-simplicial faces that continues the correspondence from
[4] between $\A$-algebras and tensor differential modules with
$\infty$-simplicial faces. Further, we show that this tensor
functor sends homotopy equivalent $\A$-al\-geb\-ras into homotopy
equivalent differential modules with $\infty$-simplicial faces.

We proceed to precise definitions and statements. All modules and
maps of modules considered in this paper are assumed to be
$K$-modules and $K$-linear maps of modules, respectively, where
$K$ is an arbitrary commutative ring with unity. \vspace{0.5cm}

\centerline{\bf \S\,1. Necessary information about}
\centerline{\bf differential modules with $\infty$-simplicial
faces.} \vspace{0.5cm}

In what follows, by a bigraded module we mean any bigraded module
$X=\{X_{n,\,m}\}$, $n\geqslant 0$, $m\geqslant 0$, and by a
differential bigraded module, or, briefly, a differential module
$(X,d)$, we mean any bigraded module $X$ endowed with a
differential $d:X_{*,\bu}\to X_{*,\bu-1}$ of bidegree $(0,-1)$.

Recall that a differential module with simplicial faces is defined
as a differential module $(X,d)$ together with a family of module
maps $\p_i:X_{n,\bu}\to X_{n-1,\bu}$, $0\leqslant i\leqslant n$,
which are maps of differential modules and satisfy the simplicial
commutation relations $\p_i\p_j=\p_{j-1}\p_i$, $i<j$. The maps
$\p_i:X_{n,\bu}\to X_{n-1,\bu}$ are called the simplicial face
operators or, more briefly, the simplicial faces of the
differential module $(X,d)$.

Now, we recall the notion of a differential module with
$\infty$-simplicial faces \cite{Lap2} (see also
\cite{Lap3}-\cite{Lap7}), which is a homotopy invariant analogue
of the notion of a differential module with simplicial faces.

Let $\Sigma_k$ be the symmetric group of permutations on a
$k$-element set. Given an arbitrary permutation
$\sigma\in\Sigma_k$ and any $k$-tuple of nonnegative integers
$(i_1,\dots,i_k)$, where $i_1<\dots<i_k$, we consider the
$k$-tuple $(\sigma(i_1),\dots,\sigma(i_k))$, where $\sigma$ acts
on the $k$-tuple $(i_1,\dots,i_k)$ in the standard way, i.e.,
permutes its components. For the $k$-tuple
$(\sigma(i_1),\dots,\sigma(i_k))$, we define a $k$-tuple
$(\widehat{\sigma(i_1)},\dots,\widehat{\sigma(i_k)})$ by the
following formulae $$\widehat{\sigma
(i_s)}=\sigma(i_s)-\alpha(\sigma(i_s)),\quad 1\leqslant s\leqslant
k,$$ where each $\alpha(\sigma(i_s))$ is the number of those
elements of $(\sigma(i_1),\dots,\sigma(i_s),\dots\sigma(i_k))$ on
the right of $\sigma(i_s)$ that are smaller than $\sigma(i_s)$.

A differential module with $\infty$-simplicial faces
$(X,d,\widetilde{\p})$ is defined as a differential module $(X,d)$
together with a family of module maps
$$\widetilde{\p}=\{\p_{(i_1,\dots ,i_k)}:X_{n,\bu}\to
X_{n-k,\bu+k-1}\},\quad 1\leqslant k\leqslant n,$$
$$i_1,\dots,i_k\in\mathbb{Z},\quad 0\leqslant
i_1<\dots<i_k\leqslant n,$$ which satisfy the relations
$$d(\p_{(i_1,\dots,i_k)})=\sum_{\sigma\in\Sigma_k}\sum_{I_{\sigma}}
(-1)^{{\rm sign}(\sigma)+1}
\p_{(\widehat{\sigma(i_1)},\dots,\widehat{\sigma(i_m)})}\,
\p_{(\widehat{\sigma(i_{m+1})},\dots,\widehat{\sigma(i_k)})},\eqno(1.1)$$
where $I_\sigma$ is the set of all partitions of the $k$-tuple
$(\widehat{\sigma(i_1)},\dots,\widehat{\sigma(i_k)})$ into two
tuples $(\widehat{\sigma(i_1)},\dots,\widehat{\sigma(i_m)})$ and
$(\widehat{\sigma(i_{m+1})},\dots,\widehat{\sigma(i_k)})$,
$1\leqslant m\leqslant k-1$, such that the conditions
$\widehat{\sigma(i_1)}<\dots<\widehat{\sigma(i_m)}$ and
$\widehat{\sigma(i_{m+1})}<\dots<\widehat{\sigma(i_k)}$ holds. The
maps $\p_{(i_1,\dots ,i_k)}\in\widetilde{\p}$ are called the
$\infty$-simplicial faces of the differential module with
$\infty$-simplicial faces $(X,d,\widetilde{\p})$.

For example, at $k=1,2,3$ the relations $(1.1)$ take,
respectively, the following view $$d(\p_{(i)})=0,\quad i\geqslant
0,\quad d(\p_{(i,j)})=\p_{(j-1)}\p_{(i)}-\p_{(i)}\p_{(j)},\quad
i<j,$$
$$d(\p_{(i_1,i_2,i_3)})=-\p_{(i_1)}\p_{(i_2,i_3)}-\p_{(i_1,i_2)}\p_{(i_3)}-
\p_{(i_3-2)}\p_{(i_1,i_2)}-$$
$$-\,\p_{(i_2-1,i_3-1)}\p_{(i_1)}+\p_{(i_2-1)}\p_{(i_1,i_3)}+\p_{(i_1,i_3-1)}\p_{(i_2)},
\quad i_1<i_2<i_3.$$

Simplest examples of differential modules with $\infty$-simplicial
faces are differential modules with simplicial faces. Indeed,
given any differential module with simplicial faces $(X,d,\p_i)$,
we can define $\widetilde{\p}=\{\p_{(i_1,\dots ,i_k)}\}:X\to X$ by
setting $\p_{(i)}=\p_i$, $i\geqslant 0$, and
$\p_{(i_1,\dots,i_k)}=0$, $k>1$, thus obtaining the differential
module with $\infty$-simplicial faces $(X,d,\widetilde{\p})$.

It is worth mentioning that the notion of an differential module
with $\infty$-simplicial faces specified above is a part of the
general notion of a differential $\infty$-simplicial module
introduced in \cite{Lap3} by using the homotopy technique of
differential Lie modules over curved colored coalgebras.

Now, we recall that a map $f:(X,d,\p_i)\to (Y,d,\p_i)$ of
differential modules with simplicial faces is defined as a map of
differential modules $f:(X,d)\to(Y,d)$ that satisfies the
relations $\p_if=f\p_i$, $i\geqslant 0$.

Let us consider the notion of a morphism of differential modules
with $\infty$-simplicial faces \cite{Lap2}, which homotopically
generalizes the notion of a map differential modules with
simplicial faces.

A morphism of differential modules with $\infty$-simplicial faces
$\widetilde{f}:(X,d,\widetilde{\p})\to (Y,d,\widetilde{\p})$ is
defined as a family of module maps
$$\widetilde{f}=\{f_{(i_1,\dots,i_k)}:X_{n,\bu}\to
Y_{n-k,\bu+k}\},\quad 0\leqslant k\leqslant n,$$
$$i_1,\dots,i_k\in\mathbb{Z},\quad 0\leqslant
i_1<\dots<i_k\leqslant n,$$ (at $k=0$ we will use the denotation
$f_{(\,\,)}$), which satisfy the relations
$$d(f_{(i_1,\dots,i_k)})=-\p_{(i_1,\dots,i_k)}f_{(\,\,)}+f_{(\,\,)}\p_{(i_1,\dots,i_k)}\,+$$
$$+\sum_{\sigma\in\Sigma_k}\sum_{I_{\sigma}}(-1)^{{\rm
sign}(\sigma)+1}\p_{(\widehat{\sigma(i_1)},\dots,
\widehat{\sigma(i_m)})}f_{(\widehat{\sigma(i_{m+1})},\dots,
\widehat{\sigma(i_k)})}\,-$$ $$-\,f_{(\widehat{\sigma(i_1)},\dots,
\widehat{\sigma(i_m)})}\p_{(\widehat{\sigma(i_{m+1})},\dots,
\widehat{\sigma(i_k)})},\eqno(1.2)$$ where $I_{\sigma}$ is the
same as in $(1.1)$. The maps $f_{(i_1,\dots
,i_k)}\in\widetilde{f}$ are called the components of the morphism
$\widetilde{f}:(X,d,\widetilde{\p})\to (Y,d,\widetilde{\p})$.

For example, at $k=0,1,2,3$ the relations $(1.2)$ take,
respectively, the following view $$d(f_{(\,\,)})=0,\qquad
d(f_{(i)})=f_{(\,\,)}\p_{(i)}-\p_{(i)}f_{(\,\,)},\quad i\geqslant
0,$$
$$d(f_{(i,j)})=-\p_{(i,j)}f_{(\,\,)}+f_{(\,\,)}\p_{(i,j)}-\p_{(i)}f_{(j)}+
\p_{(j-1)}f_{(i)}+ f_{(i)}\p_{(j)}-f_{(j-1)}\p_{(i)},\quad i<j,$$
$$d(f_{(i_1,i_2,i_3)})=-\p_{(i_1,i_2,i_3)}f_{(\,\,)}+f_{(\,\,)}\p_{(i_1,i_2,i_3)}-
\p_{(i_1)}f_{(i_2,i_3)}-\p_{(i_1,i_2)}f_{(i_3)}-
\p_{(i_3-2)}f_{(i_1,i_2)}-$$
$$-\,\p_{(i_2-1,i_3-1)}f_{(i_1)}+\p_{(i_2-1)}f_{(i_1,i_3)}+\p_{(i_1,i_3-1)}f_{(i_2)}+
f_{(i_1)}\p_{(i_2,i_3)}+f_{(i_1,i_2)}\p_{(i_3)}+$$
$$+f_{(i_3-2)}\p_{(i_1,i_2)}+
f_{(i_2-1,i_3-1)}\p_{(i_1)}-f_{(i_2-1)}\p_{(i_1,i_3)}-f_{(i_1,i_3-1)}\p_{(i_2)},\quad
i_1<i_2<i_3.$$

Under a composition $\widetilde{g}\widetilde{f}$ of morphisms of
differential modules with $\infty$-simplicial faces
$\widetilde{f}:(X,d,\widetilde{\p})\to (Y,d,\widetilde{\p})$ and
$\widetilde{g}:(Y,d,\widetilde{\p})\to (Z,d,\widetilde{\p})$ we
mean \cite{Lap2} a morphism differential modules with
$\infty$-sim\-pli\-cial faces
$\widetilde{g}\widetilde{f}=\{(gf)_{(i_1,\dots,i_k)}\}:(X,d,\widetilde{\p})\to
(Z,d,\widetilde{\p})$ whose components are defined by
$$(gf)_{(i_1,\dots,i_k)}=\sum_{\sigma\in\Sigma_k}\sum_{I'_{\sigma}}
(-1)^{{\rm sign}(\sigma)}
g_{(\widehat{\sigma(i_1)},\dots,\widehat{\sigma(i_m)})}
f_{(\widehat{\sigma(i_{m+1})},\dots,\widehat{\sigma(i_k)})},\eqno(1.3)$$
where $I'_\sigma$ is the set of all partitions of the $k$-tuple
$(\widehat{\sigma(i_1)},\dots,\widehat{\sigma(i_k)})$ into two
tuples $(\widehat{\sigma(i_1)},\dots,\widehat{\sigma(i_m)})$ and
$(\widehat{\sigma(i_{m+1})},\dots,\widehat{\sigma(i_k)})$,
$0\leqslant m\leqslant k$, such that the conditions
$\widehat{\sigma(i_1)}<\dots<\widehat{\sigma(i_m)}$ and
$\widehat{\sigma(i_{m+1})}<\dots<\widehat{\sigma(i_k)}$ holds.

For example, at $k=0,1,2,3$ the formulae $(1.3)$ take,
respectively, the following form
$$(gf)_{(\,\,)}=g_{(\,\,)}f_{(\,\,)},\qquad(gf)_{(i)}=g_{(\,\,)}f_{(i)}+g_{(i)}f_{(\,\,)},$$
$$(gf)_{(i_1,i_2)}=g_{(\,\,)}f_{(i_1,i_2)}+g_{(i_1,i_2)}f_{(\,\,)}+
g_{(i_1)}f_{(i_2)}-g_{(i_2-1)}f_{(i_1)},\qquad i_1<i_2,$$
$$(gf)_{(i_1,i_2,i_3)}=g_{(\,\,)}f_{(i_1,i_2,i_3)}+g_{(i_1,i_2,i_3)}f_{(\,\,)}+
g_{(i_1)}f_{(i_2,i_3)}+g_{(i_1,i_2)}f_{(i_3)}+$$
$$+\,g_{(i_3-2)}f_{(i_1,i_2)}+
g_{(i_2-1,i_3-1)}f_{(i_1)}-g_{(i_2-1)}f_{(i_1,i_3)}-g_{(i_1,i_3-1)}f_{(i_2)},\quad
i_1<i_2<i_3.$$

Given any differential module with $\infty$-simplicial faces
$(X,d,\widetilde{\p})$, there is an identity morphism
$$\widetilde{1}_X=\{(1_X)_{(i_1,\dots,i_k)}\}:(X,d,\widetilde{\p})\to
(X,d,\widetilde{\p}),$$ where $(1_X)_{(\,\,)}=1_X$ is the identity
map of the module $X$ and $(1_X)_{(i_1,\dots,i_k)}=0$ is the zero
map of modules for each $k\geqslant 1$. Thus, the class of all
differential modules with $\infty$-simplicial faces and their
morphisms is a category.

Now, we recall that a homotopy between morphisms
$f,g:(X,d,\p_i)\to(Y,d,\p_i)$ of differential modules with
simplicial faces is defines as a homotopy $h:X_{*,\bu}\to
Y_{*,\bu+1}$ between morphisms of differential modules
$f,g:(X,d)\to (Y,d)$ satisfies the relations $$\p_ih+h\p_i=0,\quad
i\geqslant 0.$$

Let us consider the notion of a homotopy between morphisms of
differential modules with $\infty$-simplicial faces \cite{Lap2},
which homotopically generalizes the notion of a homotopy between
morphisms of differential modules with simplicial faces.

A homotopy $\widetilde{h}:(X,d,\widetilde{\p})\to
(Y,d,\widetilde{\p})$ between morphisms of differential modules
with $\infty$-simplicial faces
$\widetilde{f},\widetilde{g}:(X,d,\widetilde{\p})\to
(Y,d,\widetilde{\p})$ is defined as a family of module maps
$$\widetilde{h}=\{h_{(i_1,\dots,i_k)}:X_{n,\bu}\to
Y_{n-k,\bu+k+1}\},\quad 0\leqslant k\leqslant n,$$
$$i_1,\dots,i_k\in\mathbb{Z},\quad 0\leqslant
i_1<\dots<i_k\leqslant n,$$ (at $k=0$ we will use the denotation
$h_{(\,\,)}$), which satisfy the relations
$$d(h_{(i_1,\dots,i_k)})=f_{(i_1,\dots,i_k)}-g_{(i_1,\dots,i_k)}-
\partial_{(i_1,\dots,i_k)}h_{(\,\,)}-h_{(\,\,)}\partial_{(i_1,\dots,i_k)}\,+$$
$$+\sum_{\sigma\in\Sigma_k}\sum_{I_{\sigma}}(-1)^{{\rm
sign}(\sigma)+1}
\partial_{(\widehat{\sigma(i_1)},\dots,\widehat{\sigma(i_m)})}
h_{(\widehat{\sigma(i_{m+1})},\dots,\widehat{\sigma(i_k)})}\,+$$
$$+\,h_{(\widehat{\sigma(i_1)},\dots,\widehat{\sigma(i_m)})}
\partial_{(\widehat{\sigma(i_{m+1})},\dots,\widehat{\sigma(i_k)})},\eqno(1.4)$$
where $I_{\sigma}$ is the same as in $(1.1)$. The maps
$h_{(i_1,\dots ,i_k)}\in\widetilde{h}$ are called the components
of the homotopy $\widetilde{h}:(X,d,\widetilde{\p})\to
(Y,d,\widetilde{\p})$.

For example, at $k=0,1,2,3$ the relations $(1.4)$ take,
respectively, the following view
$$d(h_{(\,\,)})=f_{(\,\,)}-g_{(\,\,)},\quad
d(h_{(i)})=f_{(i)}-g_{(i)}-\p_{(i)}h_{(\,\,)}-h_{(\,\,)}\p_{(i)},\quad
i\geqslant 0,$$ $$d(h_{(i,j)})=f_{(i,j)}-g_{(i,j)}-
\p_{(i,j)}h_{(\,\,)}-h_{(\,\,)}\p_{(i,j)}-\p_{(i)}h_{(j)}+\p_{(j-1)}h_{(i)}-
h_{(i)}\p_{(j)}+h_{(j-1)}\p_{(i)},~i<j,$$
$$d(h_{(i_1,i_2,i_3)})=f_{(i_1,i_2,i_3)}-g_{(i_1,i_2,i_3)}-\p_{(i_1,i_2,i_3)}h_{(\,\,)}-h_{(\,\,)}\p_{(i_1,i_2,i_3)}-
\p_{(i_1)}f_{(i_2,i_3)}-\p_{(i_1,i_2)}f_{(i_3)}\,-$$
$$-\,\p_{(i_3-2)}f_{(i_1,i_2)}-\p_{(i_2-1,i_3-1)}f_{(i_1)}+\p_{(i_2-1)}f_{(i_1,i_3)}+\p_{(i_1,i_3-1)}f_{(i_2)}-
h_{(i_1)}\p_{(i_2,i_3)}-h_{(i_1,i_2)}\p_{(i_3)}\,-$$
$$-\,h_{(i_3-2)}\p_{(i_1,i_2)}-
h_{(i_2-1,i_3-1)}\p_{(i_1)}+h_{(i_2-1)}\p_{(i_1,i_3)}+h_{(i_1,i_3-1)}\p_{(i_2)},\quad
i_1<i_2<i_3.$$

By using the notion of a homotopy between morphisms of
differential modules with $\infty$-simplicial faces the notion
homotopy equivalent differential modules with
$\infty$-sim\-plicial faces is introduced in the usual way.

In conclusion of this paragraph it is worth mentioning that the
homotopy invariance of the structure of a differential module with
$\infty$-simplicial faces under homotopy equivalences of the type
of an SDR-data of differential modules was established in
\cite{Lap2}.
\vspace{0.5cm}

\centerline{\bf \S\,2. Necessary information about $\A$-algebras.}
\vspace{0.5cm}

Recall, following \cite{Kad} (see also \cite{S}), that an
$\A$-algebra $(A,d,\pi_n)$ is any differential module $(A,d)$ with
$A=\{A_n\}$, $n\in\mathbb{Z}$, $n\geqslant 0$, $d:A_\bu\to
A_{\bu-1}$, equipped with a family of maps $\{\pi_n:(A^{\otimes
(n+2)})_\bu\to A_{\bu+n}\}$, $n\geqslant 0$, satisfying the
following relations for any integer $n\geqslant -1$:
$$d(\pi_{n+1})=\sum\limits_{m=0}^n\sum_{t=1}^{m+2}(-1)^{t(n-m+1)+n+1}
\pi_m(\underbrace{1\otimes\dots\otimes
1}_{t-1}\otimes\,\pi_{n-m}\otimes\underbrace{1\otimes\dots \otimes
1}_{m-t+2}),\eqno(2.1)$$ where
$d(\pi_{n+1})=d\pi_{n+1}+(-1)^n\pi_{n+1}d$.

For example, at $n=-1,0,1$ the relations $(2.1)$ take the forms
$$d(\pi_0)=0,\quad d(\pi_1)=\pi_0(\pi_0\otimes
1)-\pi_0(1\otimes\pi_0),$$ $$d(\pi_2)=\pi_0(\pi_1\otimes
1)+\pi_0(1\otimes\pi_1)-\pi_1(\pi_0\otimes 1\otimes
1)+\pi_1(1\otimes\pi_0\otimes 1)-\pi_1(1\otimes 1\otimes\pi_0).$$

For greater clarity and assurance that in $(2.1)$ the sign
$(-1)^{t(n-m+1)+n+1}$ is written correctly, let us describe the
procedure of finding of this sign.

Given any differential module $(A,d)$ with $A=\{A_n\}$,
$n\in\mathbb{Z}$, $n\geqslant 0$, $d:A_\bu\to A_{\bu-1}$, consider
the suspension $(SA,d)$ over $(A,d)$, where $(SA)_{\bu+1}=A_\bu$
and the differential $d:(SA)_\bu\to (SA)_{\bu-1}$ at any element
$[a]\in (SA)_\bu$ is defined by $d([a])=[d(a)]$. It is well known
(see, for example, \cite{Smir}) that introduction of the structure
of an $\A$-algebra $(A,d,\pi_n)$ on the differential module
$(A,d)$ is equivalent to consideration of a family of module maps
$$\{\pi(n):(SA)^{\otimes (n+2)}_\bu\to (SA)_{\bu-1}\},\quad
n\geqslant 0,$$ satisfying the following relations for any integer
$n\geqslant -1$: $$d(\pi(n+1))+\sum\limits_{m=0}^n\sum_{t=1}^{m+2}
\pi(m)(\underbrace{1\otimes\dots\otimes
1}_{t-1}\otimes\,\pi(n-m)\otimes\underbrace{1\otimes\dots \otimes
1}_{m-t+2})=0\,,$$ where $d(\pi(n+1))=d\pi(n+1)+\pi(n+1)d$.
Indeed, if we consider the maps of differential modules
$\eta:((SA)_\bu,d)\to (A_{\bu-1},d)$ and $\xi:(A_\bu,d)\to
((SA)_{\bu+1},d)$ defined by $\eta([a])=a$ and $\xi(a)=[a]$, then
the specified above maps $\pi_n$ and $\pi(n)$ define each other by
$\pi(n)=\xi\pi_n\eta^{\otimes (n+2)}$. By using equalities
$d\eta=\eta d$, $\deg(\eta)=-1$, $\deg(\pi(n))=-1$,
$\deg(\pi_n)=n$, $\eta\pi(n)=\pi_n\eta^{\otimes (n+2)}$, and also
by applying the Koszul rule of computing of signs, we obtain
$$\eta\left(d(\pi(n+1))+\sum\limits_{m=0}^n\sum_{t=1}^{m+2}
\pi(m)(\underbrace{1\otimes\dots\otimes
1}_{t-1}\otimes\,\pi(n-m)\otimes\underbrace{1\otimes\dots \otimes
1}_{m-t+2})\right)=$$ $$=d(\pi_{n+1})\eta^{\otimes
(n+3)}+\sum\limits_{m=0}^n\sum_{t=1}^{m+2} \pi_m\eta^{\otimes
(m+2)}(\underbrace{1\otimes\dots\otimes
1}_{t-1}\otimes\,\pi(n-m)\otimes\underbrace{1\otimes\dots \otimes
1}_{m-t+2})=$$ $$=d(\pi_{n+1})\eta^{\otimes
(n+3)}+\sum\limits_{m=0}^n\sum_{t=1}^{m+2}(-1)^{m-t+2}
\pi_m(\eta^{\otimes (t-1)}\otimes\eta\pi(n-m)\otimes\eta^{\otimes
(m-t+2)})=$$ $$=d(\pi_{n+1})\eta^{\otimes
(n+3)}+\sum\limits_{m=0}^n\sum_{t=1}^{m+2}(-1)^{m-t}
\pi_m(\eta^{\otimes (t-1)}\otimes \pi_{n-m}\eta^{\otimes (n-m+2)}
\otimes\eta^{\otimes (m-t+2)})=$$
$$=\left(d(\pi_{n+1})+\sum\limits_{m=0}^n\sum_{t=1}^{m+2}(-1)^\varepsilon
\pi_m(1^{\otimes (t-1)}\otimes\pi_{n-m} \otimes 1^{\otimes
(m-t+2)})\right)\eta^{\otimes (n+3)}=0,$$ where
$\varepsilon=m-t+(t-1)(n-m)$. It follows that the relations
$(2.1)$ holds because the map $\eta$ is a isomorphism of degree
$(-1)$ of graded modules and the congruence $m-t+(t-1)(n-m)\equiv
t(n-m-1)+n\,{\rm mod}(2)$ is true.

Now, recall \cite{Kad} that a morphism of $\A$-algebras
$f:(A,d,\pi_n)\to (A',d,\pi_n)$ is defined as a family of module
maps $f=\{f_n:(A^{\otimes(n+1)})_\bu\to
A'_{\bu+n}~|~n\in\mathbb{Z},~n\geqslant 0\}$, which, for all
integers $n\geqslant -1$, satisfy the relations
$$d(f_{n+1})=\sum_{m=0}^n\sum_{t=1}^{m+1}(-1)^{t(n-m+1)+n+1}f_m(\underbrace{1\otimes\dots\otimes
1}_{t-1}\otimes\,\pi_{n-m}\otimes\underbrace{1\otimes\dots\otimes
1}_{m-t+1})\,-$$
$$-\sum_{m=0}^n\sum_{\,\,\,J_{n-m}}(-1)^{\varepsilon}\pi_m(f_{n_1}\otimes
f_{n_2}\otimes\dots\otimes f_{n_{m+2}}),\eqno(2.2)$$ where
$d(f_{n+1})=df_{n+1}+(-1)^nf_{n+1}d$ and $$J_{n-m}=\{n_1\geqslant
0,n_2\geqslant 0,\dots,n_{m+2}\geqslant
0~|~n_1+n_2+\dots+n_{m+2}=n-m\},$$
$$\varepsilon=\sum_{i=1}^{m+1}(n_i+1)(n_{i+1}+\dots+n_{m+2}).$$

For example, at $n=-1,0,1$ the relations $(2.2)$ take,
respectively, the following view $$d(f_0)=0,\quad
d(f_1)=f_0\pi_0-\pi_0(f_0\otimes f_0),$$
$$d(f_2)=f_0\pi_1-f_1(\pi_0\otimes
1)+f_1(1\otimes\pi_0)-\pi_0(f_1\otimes f_0)+\pi_0(f_0\otimes
f_1)-\pi_1(f_0\otimes f_0\otimes f_0).$$

For greater clarity and assurance that in $(2.2)$ the signs
$(-1)^{t(n-m+1)+n+1}$ and $(-1)^{\varepsilon}$ are written
correctly, let us describe the procedure of finding of these
signs.

As above, given any differential module $(A,d)$ with $A=\{A_n\}$,
$n\in\mathbb{Z}$, $n\geqslant 0$, $d:A_\bu\to A_{\bu-1}$, consider
the suspension $(SA,d)$ over $(A,d)$ and the maps of differential
modules $\eta:(SA)_\bu\to A_{\bu-1}$ and $\xi:A_\bu\to
(SA)_{\bu+1}$. It is well known (see, for example, \cite{Smir})
that consideration of the morphism $f=\{f_n\}:(A,d,\pi_n)\to
(A',d,\pi_n)$ of $\A$-al\-geb\-ras is equivalent to consideration
of the family of module maps $$\{f(n):(SA)^{\otimes (n+1)}_\bu\to
(SA')_\bu\},\quad n\geqslant 0,$$ satisfying the following
relations for any integer $n\geqslant -1$:
$$d(f(n+1))=\sum_{m=0}^n\sum_{t=1}^{m+1}f(m)(\underbrace{1\otimes\dots\otimes
1}_{t-1}\otimes\,\pi(n-m)\otimes\underbrace{1\otimes\dots\otimes
1}_{m-t+1})\,-$$
$$-\sum_{m=0}^n\sum_{\,\,\,J_{n-m}}\pi(m)(f(n_1)\otimes
f(n_2)\otimes\dots\otimes f(n_{m+2})),$$ where $J_{n-m}$ is the
same as in $(2.2)$ and $d(f(n+1))=df(n+1)-f(n+1)d$. The maps
$\pi(n)$, $n\geqslant 0$, were defined above. The specified above
maps $f_n$ and $f(n)$ define each other by $f(n)=\xi
f_n\eta^{\otimes (n+1)}$, $n\geqslant 0$. Since $d\eta=\eta d$ and
$\eta f(n)=f_n\eta^{\otimes (n+1)}$ for all $n\geqslant 0$, we
have $$\eta d(f(n+1))=d(f_{n+1})\eta^{\otimes (n+2)}.$$ By using
equalities $\deg(\eta)=-1$, $\deg(\pi(n))=-1$, $\deg(\pi_n)=n$,
$\eta f(n)=f_n\eta^{\otimes (n+1)}$,
$\eta\pi(n)=\pi_n\eta^{\otimes (n+2)}$, and also by applying the
Koszul rule of computing of signs, we obtain
$$\eta\left(\sum_{m=0}^n\sum_{t=1}^{m+1}f(m)(\underbrace{1\otimes\dots\otimes
1}_{t-1}\otimes\,\pi(n-m)\otimes\underbrace{1\otimes\dots\otimes
1}_{m-t+1})\right)=$$
$$=\sum_{m=0}^n\sum_{t=1}^{m+1}f_m\eta^{\otimes (m+1)}
(\underbrace{1\otimes\dots\otimes
1}_{t-1}\otimes\,\pi(n-m)\otimes\underbrace{1\otimes\dots\otimes
1}_{m-t+1})=$$
$$=\sum_{m=0}^n\sum_{t=1}^{m+1}(-1)^{m-t+1}f_m(\underbrace{\eta\otimes\dots\otimes
\eta}_{t-1}\otimes\,\eta\pi(n-m)\otimes\underbrace{\eta\otimes\dots\otimes
\eta}_{m-t+1})=$$
$$=\sum_{m=0}^n\sum_{t=1}^{m+1}(-1)^{m-t+1}f_m(\underbrace{\eta\otimes\dots\otimes
\eta}_{t-1}\otimes\,\pi_{n-m}\eta^{\otimes (n-m+2)}
\otimes\underbrace{\eta\otimes\dots\otimes \eta}_{m-t+1})=$$
$$=\left(\sum_{m=0}^n\sum_{t=1}^{m+1}(-1)^{m-t+1+(t-1)(n-m)}f_m(\underbrace{1\otimes\dots\otimes
1}_{t-1}\otimes\,\pi_{n-m} \otimes\underbrace{1\otimes\dots\otimes
1}_{m-t+1})\right)\eta^{\otimes (n+2)}=$$
$$=\left(\sum_{m=0}^n\sum_{t=1}^{m+1}(-1)^{t(n-m+1)+n+1}f_m(\underbrace{1\otimes\dots\otimes
1}_{t-1}\otimes\,\pi_{n-m} \otimes\underbrace{1\otimes\dots\otimes
1}_{m-t+1})\right)\eta^{\otimes (n+2)}.$$ By using equalities
$\deg(\eta)=-1$, $\deg(f(n))=0$, $\deg(f_n)=n$,
$\eta\pi(n)=\pi_n\eta^{\otimes (n+2)}$, $\eta
f(n)=f_n\eta^{\otimes (n+1)}$, and also by applying the Koszul
rule of computing of signs, we have
$$\eta\left(-\sum_{m=0}^n\sum_{\,\,\,J_{n-m}}\pi(m)(f(n_1)\otimes
f(n_2)\otimes\dots\otimes f(n_{m+2}))\right)=$$
$$=-\sum_{m=0}^n\sum_{\,\,\,J_{n-m}}\pi_m(f_{n_1}\eta^{\otimes
(n_1+1)}\otimes f_{n_2}\eta^{\otimes (n_2+1)}\otimes\dots\otimes
f_{n_{m+2}}\eta^{\otimes (n_{m+2}+1)})=$$
$$=-\sum_{m=0}^n\sum_{\,\,\,J_{n-m}}(-1)^{(n_{m+1}+1)n_{m+2}}\pi_m(f_{n_1}\eta^{\otimes
(n_1+1)}\otimes\dots\otimes f_{n_m}\eta^{\otimes
(n_m+1)}\,\otimes$$ $$\otimes\, f_{n_{m+1}}\otimes
f_{n_{m+2}})(1^{\otimes (n_1+1)}\otimes\dots\otimes 1^{\otimes
(n_m+1)}\otimes\eta^{\otimes (n_{m+1}+1)}\otimes\eta^{\otimes
(n_{m+2}+1)})=$$
$$=-\sum_{m=0}^n\sum_{\,\,\,J_{n-m}}(-1)^{(n_{m+1}+1)n_{m+2}+(n_m+1)(n_{m+1}+n_{m+2})}\pi_m(f_{n_1}\eta^{\otimes
(n_1+1)}\otimes\dots$$ $$\dots\otimes f_{n_{m-1}}\eta^{\otimes
(n_{m-1}+1)}\otimes f_{n_m}\otimes f_{n_{m+1}}\otimes
f_{n_{m+2}})(1^{\otimes (n_1+1)}\otimes\dots$$

$$\dots\otimes 1^{\otimes (n_{m-1}+1)}\otimes\eta^{\otimes
(n_m+1)}\otimes\eta^{\otimes (n_{m+1}+1)}\otimes\eta^{\otimes
(n_{m+2}+1)})=\dots$$
$$\dots=\left(-\sum_{m=0}^n\sum_{\,\,\,J_{n-m}}(-1)^\varepsilon\pi_m(f_{n_1}\otimes
f_{n_2}\otimes\dots\otimes f_{n_{m+2}})\right)\eta^{\otimes
(n+2)},$$ where
$$\varepsilon=\sum_{i=1}^{m+1}(n_i+1)(n_{i+1}+\dots+n_{m+2}).$$
The obtained above equalities follows that the relations $(2.2)$
holds because the map $\eta$ is a isomorphism of degree $(-1)$ of
graded modules.

Under a composition of morphisms of $\A$-algebras
$f:(A,d,\pi_n)\to(A',d,\pi_n)$ and
$g:(A',d,\pi_n)\to(A'',d,\pi_n)$ we mean \cite{Kad} a morphism of
$\A$-algebras $$gf=\{(gf)_n\}:(A,d,\pi_n)\to(A'',d,\pi_n)$$
defined by $$(gf)_n=\sum_{m=0}^n\sum_{J_{n-m}}(-1)^\varepsilon
g_m(f_{n_1}\otimes f_{n_2}\otimes\dots\otimes
f_{n_{m+1}}),\eqno(2.3)$$ where $J_{n-m}=\{n_1\geqslant
0,n_2\geqslant 0,\dots,n_{m+1}\geqslant
0~|~n_1+n_2+\dots+n_{m+1}=n-m\}$ and
$$\varepsilon=\sum_{i=1}^m(n_i+1)(n_{i+1}+\dots+n_{m+1}).$$

For example, at $n=0,1,2$ the formulae $(2.3)$ take, respectively,
the following view $$(gf)_0=g_0f_0,\quad
(gf)_1=g_0f_1+g_1(f_0\otimes f_0),$$
$$(gf)_2=g_0f_2-g_1(f_0\otimes f_1)+g_1(f_1\otimes
f_0)+g_2(f_0\otimes f_0\otimes f_0).$$

For greater clarity and assurance that in $(2.3)$ the sign
$(-1)^\varepsilon$ is written correctly, let us describe the
procedure of finding of this sign.

As above, given any differential module $(A,d)$ with $A=\{A_n\}$,
$n\in\mathbb{Z}$, $n\geqslant 0$, $d:A_\bu\to A_{\bu-1}$, consider
the suspension $(SA,d)$ over $(A,d)$ and the map of differential
modules $\eta:(SA)_\bu\to A_{\bu-1}$. It is well known (see, for
example, \cite{Smir}) that consideration of the composition
$gf=\{(gf)_n\}:(A,d,\pi_n)\to (A'',d,\pi_n)$ of morphisms of
$\A$-al\-geb\-ras $f=\{f_n\}:(A,d,\pi_n)\to (A',d,\pi_n)$ and
$g=\{g_n\}:(A',d,\pi_n)\to (A'',d,\pi_n)$ is equivalent to
consideration of the composition $\{(gf)(n):(SA)^{\otimes
(n+1)}_\bu\to (SA'')_\bu\}$ of families of maps
$\{f(n):(SA)^{\otimes (n+1)}_\bu\to (SA')_\bu\}$ and
$\{g(n):(SA')^{\otimes (n+1)}_\bu\to (SA'')_\bu\}$ corresponding
to the morphisms of $\A$-algebras $f$ and $g$. The composition
$\{(gf)(n)\}$ of families of maps $\{f(n)\}$ and $\{g(n)\}$ is
defined by $$(gf)(n)=\sum_{m=0}^n\sum_{J_{n-m}}g(m)(f(n_1)\otimes
f(n_2)\otimes\dots\otimes f(n_{m+1})),\quad n\geqslant 0,$$ where
$J_{n-m}$ is the same as in $(2.3)$. By using equalities
$$\deg(\eta)=-1,\quad \deg(f(n))=\deg(g(n))=0,\quad
\deg(f_n)=\deg(g_n)=n,$$ $$\eta f(n)=f_n\eta^{\otimes (n+1)},\quad
\eta g(n)=g_n\eta^{\otimes (n+1)},\quad \eta
(gf)(n)=(gf)_n\eta^{\otimes (n+1)},$$ and also by applying the
Koszul rule of computing of signs, we obtain
$$\eta\left((gf)(n)-\sum_{m=0}^n\sum_{J_{n-m}}g(m)(f(n_1)\otimes
f(n_2)\otimes\dots\otimes f(n_{m+1}))\right)=$$
$$=(gf)_n\eta^{\otimes
(n+1)}-\sum_{m=0}^n\sum_{\,\,\,J_{n-m}}g_m(f_{n_1}\eta^{\otimes
(n_1+1)}\otimes f_{n_2}\eta^{\otimes (n_2+1)}\otimes\dots\otimes
f_{n_{m+1}}\eta^{\otimes (n_{m+1}+1)})=$$ $$=(gf)_n\eta^{\otimes
(n+1)}-\sum_{m=0}^n\sum_{\,\,\,J_{n-m}}(-1)^{(n_m+1)n_{m+1}}g_m(f_{n_1}\eta^{\otimes
(n_1+1)}\otimes\dots\otimes f_{n_{m-1}}\eta^{\otimes
(n_{m-1}+1)}\,\otimes$$ $$\otimes\, f_{n_{m}}\otimes
f_{n_{m+1}})(1^{\otimes (n_1+1)}\otimes\dots\otimes 1^{\otimes
(n_{m-1}+1)}\otimes\eta^{\otimes (n_{m}+1)}\otimes\eta^{\otimes
(n_{m+1}+1)})=$$ $$=(gf)_n\eta^{\otimes
(n+1)}-\sum_{m=0}^n\sum_{\,\,\,J_{n-m}}(-1)^{(n_m+1)n_{m+1}+(n_{m-1}+1)(n_m+n_{m+1})}g_m(f_{n_1}\eta^{\otimes
(n_1+1)}\otimes\dots$$ $$\dots\otimes f_{n_{m-2}}\eta^{\otimes
(n_{m-2}+1)}\otimes f_{n_{m-1}}\otimes f_{n_m}\otimes
f_{n_{m+1}})(1^{\otimes (n_1+1)}\otimes\dots$$

$$\dots\otimes 1^{\otimes (n_{m-2}+1)}\otimes\eta^{\otimes
(n_{m-1}+1)}\otimes\eta^{\otimes (n_m+1)}\otimes\eta^{\otimes
(n_{m+1}+1)})=\dots$$
$$\dots=\left((gf)_n-\sum_{m=0}^n\sum_{\,\,\,J_{n-m}}(-1)^\varepsilon
g_m(f_{n_1}\otimes f_{n_2}\otimes\dots\otimes
f_{n_{m+1}})\right)\eta^{\otimes (n+1)}=0,$$ where
$$\varepsilon=\sum_{i=1}^m(n_i+1)(n_{i+1}+\dots+n_{m+1}).$$ The
obtained above equalities follows that the relations $(2.3)$ holds
because the map $\eta$ is a isomorphism of degree $(-1)$ of graded
modules.

Given any $\A$-algebra $(A,d,\pi_n)$, there is an identity
morphism $$1_A=\{(1_A)_n\}:(A,d,\pi_n)\to (A,d,\pi_n),$$ where
$(1_A)_0$ is the identity map of the module $A$ and $(1_A)_n=0$ is
the zero map of modules for each $n\geqslant 1$. Thus, the class
of all $\A$-algebras and their morphisms is a category.

Now, recall (see, for example \cite{Smir}) that a homotopy
$h:(A,d,\pi_n)\to (A',d,\pi_n)$ between morphisms of $\A$-algebras
$f=\{f_n\},g=\{g_n\}:(A,d,\pi_n)\to (A',d,\pi_n)$ is defined as a
family of module maps $h=\{h_n:(A^{\otimes(n+1)})_\bu\to
A'_{\bu+n+1}~|~n\in\mathbb{Z},~n\geqslant 0\}$, which, for all
integers $n\geqslant -1$, satisfy the relations
$$d(h_{n+1})=f_{n+1}-g_{n+1}+\sum_{m=0}^n\sum_{t=1}^{m+1}(-1)^{t(n-m+1)+n}h_m(\underbrace{1\otimes\dots\otimes
1}_{t-1}\otimes\,\pi_{n-m}\otimes\underbrace{1\otimes\dots\otimes
1}_{m-t+1})\,+$$
$$+\sum_{m=0}^n\sum_{\,\,\,J_{n-m}}\sum_{i=1}^{m+2}(-1)^\varrho\pi_m(g_{n_1}\otimes\dots\otimes
g_{n_{i-1}}\otimes h_{n_i}\otimes f_{n_{i+1}}\otimes\dots\otimes
f_{n_{m+2}}),\eqno(2.4)$$ where
$d(h_{n+1})=dh_{n+1}+(-1)^{n+1}h_{n+1}d$ and
$$J_{n-m}=\{n_1\geqslant 0,n_2\geqslant 0,\dots,n_{m+2}\geqslant
0~|~n_1+n_2+\dots+n_{m+2}=n-m\},$$
$$\varrho=m+\sum_{k=1}^{m+1}(n_k+1)(n_{k+1}+\dots+n_{m+2})+\sum_{s=1}^{i-1}n_s.$$

For example, at $n=-1,0,1$ the relations $(2.4)$ take,
respectively, the following view $$d(h_0)=f_0-g_0,\quad
d(h_1)=f_1-g_1-h_0\pi_0+\pi_0(h_0\otimes f_0)+\pi_0(g_0\otimes
f_0),$$ $$d(h_2)=f_2-g_2-h_0\pi_1+h_1(\pi_0\otimes
1)-h_1(1\otimes\pi_0)-\pi_0(h_0\otimes f_1)-\pi_0(g_0\otimes
h_1)\,+$$ $$+\,\pi_0(h_1\otimes f_0)-\pi_0(g_1\otimes
h_0)-\pi_1(h_0\otimes f_0\otimes f_0)-\pi_1(g_0\otimes h_0\otimes
f_0)-\pi_1(g_0\otimes g_0\otimes h_0).$$

For greater clarity and assurance that in $(2.4)$ the signs
$(-1)^{t(n-m+1)+n}$ and $(-1)^{\varrho}$ are written correctly,
let us describe the procedure of finding of these signs.

As above, given any differential module $(A,d)$ with $A=\{A_n\}$,
$n\in\mathbb{Z}$, $n\geqslant 0$, $d:A_\bu\to A_{\bu-1}$, consider
the suspension $(SA,d)$ over $(A,d)$ and the maps of differential
modules $\eta:(SA)_\bu\to A_{\bu-1}$ and $\xi:A_\bu\to
(SA)_{\bu+1}$. It is well known (see, for example, \cite{Smir})
that consideration of the homotopy $h=\{h_n\}:(A,d,\pi_n)\to
(A',d,\pi_n)$ between morphisms of $\A$-al\-geb\-ras
$f=\{f_n\},g=\{g_n\}:(A,d,\pi_n)\to (A',d,\pi_n)$ is equivalent to
consideration of the family of module maps $$\{h(n):(SA)^{\otimes
(n+1)}_\bu\to (SA')_{\bu+1}\},\quad n\geqslant 0,$$ satisfying the
following relations for any integer $n\geqslant -1$:
$$d(h(n+1))=f(n+1)-g(n+1)\,-$$
$$-\sum_{m=0}^n\sum_{t=1}^{m+1}h(m)(\underbrace{1\otimes\dots\otimes
1}_{t-1}\otimes\,\pi(n-m)\otimes\underbrace{1\otimes\dots\otimes
1}_{m-t+1})\,-$$
$$-\sum_{m=0}^n\sum_{\,\,\,J_{n-m}}\sum_{i=1}^{m+2}\pi(m)(g(n_1)\otimes\dots\otimes
g(n_{i-1})\otimes h(n_i)\otimes f(n_{i+1})\otimes\dots\otimes
f(n_{m+2})),$$ where $J_{n-m}$ is the same as in $(2.4)$ and
$d(h(n+1))=dh(n+1)+h(n+1)d$. The maps $\pi(n)$, $f(n)$, $g(n)$,
$n\geqslant 0$, were considered above. The specified above maps
$h_n$ and $h(n)$ define each other by $h(n)=\xi h_n\eta^{\otimes
(n+1)}$, $n\geqslant 0$. Since $d\eta=\eta d$ and $\eta
h(n)=h_n\eta^{\otimes (n+1)}$ for all $n\geqslant 0$, we have
$$\eta d(h(n+1))=d(h_{n+1})\eta^{\otimes (n+2)}.$$ It is obvious
that the equality
$$\eta(f(n+1)-g(n+1))=(f_{n+1}-g_{n+1})\eta^{\otimes (n+2)}$$
holds. By using $\deg(\eta)=-1$, $\deg(\pi(n))=-1$,
$\deg(\pi_n)=n$, $\eta h(n)=h_n\eta^{\otimes (n+1)}$,
$\eta\pi(n)=\pi_n\eta^{\otimes (n+2)}$, and also by applying the
Koszul rule of computing of signs, we obtain
$$\eta\left(-\sum_{m=0}^n\sum_{t=1}^{m+1}h(m)(\underbrace{1\otimes\dots\otimes
1}_{t-1}\otimes\,\pi(n-m)\otimes\underbrace{1\otimes\dots\otimes
1}_{m-t+1})\right)=$$
$$=-\sum_{m=0}^n\sum_{t=1}^{m+1}h_m\eta^{\otimes (m+1)}
(\underbrace{1\otimes\dots\otimes
1}_{t-1}\otimes\,\pi(n-m)\otimes\underbrace{1\otimes\dots\otimes
1}_{m-t+1})=$$
$$=\sum_{m=0}^n\sum_{t=1}^{m+1}(-1)^{m-t}h_m(\underbrace{\eta\otimes\dots\otimes
\eta}_{t-1}\otimes\,\eta\pi(n-m)\otimes\underbrace{\eta\otimes\dots\otimes
\eta}_{m-t+1})=$$
$$=\sum_{m=0}^n\sum_{t=1}^{m+1}(-1)^{m-t}h_m(\underbrace{\eta\otimes\dots\otimes
\eta}_{t-1}\otimes\,\pi_{n-m}\eta^{\otimes (n-m+2)}
\otimes\underbrace{\eta\otimes\dots\otimes \eta}_{m-t+1})=$$
$$=\left(\sum_{m=0}^n\sum_{t=1}^{m+1}(-1)^{m-t+(t-1)(n-m)}h_m(\underbrace{1\otimes\dots\otimes
1}_{t-1}\otimes\,\pi_{n-m} \otimes\underbrace{1\otimes\dots\otimes
1}_{m-t+1})\right)\eta^{\otimes (n+2)}=$$
$$=\left(\sum_{m=0}^n\sum_{t=1}^{m+1}(-1)^{t(n-m+1)+n}h_m(\underbrace{1\otimes\dots\otimes
1}_{t-1}\otimes\,\pi_{n-m} \otimes\underbrace{1\otimes\dots\otimes
1}_{m-t+1})\right)\eta^{\otimes (n+2)}.$$ By using equalities
$$\deg(\eta)=-1,\quad \deg(h(n))=1,\quad \deg(h_n)=n+1,\quad
\deg(f(n))=\deg(g(n))=0,$$ $$\deg(f_n)=\deg(g_n)=n,\quad \eta
h(n)=h_n\eta^{\otimes (n+1)},\quad \eta\pi(n)=\pi_n\eta^{\otimes
(n+2)},$$ $$\eta f(n)=f_n\eta^{\otimes (n+1)},\quad \eta
g(n)=g_n\eta^{\otimes (n+1)},$$ and also by applying the Koszul
rule of computing of signs, we have
$$\eta\!\left(\!-\!\sum_{m=0}^n\sum_{\,\,\,J_{n-m}}\sum_{i=1}^{m+2}\pi(m)(g(n_1)\otimes\dots\otimes
g(n_{i-1})\otimes h(n_i)\otimes f(n_{i+1})\otimes\dots\otimes
f(n_{m+2}))\!\right)\!\!=$$
$$=-\sum_{m=0}^n\sum_{\,\,\,J_{n-m}}\sum_{i=1}^{m+2}(-1)^{m-i+2}\pi_m(\underbrace{\eta
g(n_1)\otimes\dots\otimes \eta g(n_{i-1})}_{i-1}\otimes\,\eta
h(n_i)\,\otimes$$ $$\otimes\underbrace{\eta
f(n_{i+1})\otimes\dots\otimes \eta f(n_{m+2})}_{m-i+2})=$$
$$=\sum_{m=0}^n\sum_{\,\,\,J_{n-m}}\sum_{i=1}^{m+2}(-1)^{m-i+1}\pi_m(g_{n_1}\eta^{\otimes
(n_1+1)}\otimes\dots\otimes g_{n_{i-1}}\eta^{\otimes (n_{i-1}+1)}
\otimes h_{n_i}\eta^{\otimes (n_i+1)}\,\otimes$$ $$\otimes\,
f_{n_{i+1}}\eta^{\otimes (n_{i+1}+1)}\otimes\dots\otimes
f_{n_{m+2}}\eta^{\otimes (n_{m+2}+1)})=$$
$$=\sum_{m=0}^n\sum_{\,\,\,J_{n-m}}\sum_{i=1}^{m+2}(-1)^{m-i+1+\alpha}\pi_m(g_{n_1}\eta^{\otimes
(n_1+1)}\otimes\dots\otimes g_{n_{i-1}}\eta^{\otimes (n_{i-1}+1)}
\otimes h_{n_i}\otimes f_{n_{i+1}}\otimes\dots$$ $$\dots\otimes
f_{n_{m+2}})(1^{\otimes (n_1+1)}\otimes\dots\otimes 1^{\otimes
(n_{i-1}+1)}\otimes \eta^{\otimes (n_i+1)}\otimes \eta^{\otimes
(n_{i+1}+1)}\otimes\dots\otimes \eta^{\otimes (n_{m+2}+1)})=$$
$$=\sum_{m=0}^n\sum_{\,\,\,J_{n-m}}\sum_{i=1}^{m+2}(-1)^{m-i+1+\alpha+\beta}\pi_m(g_{n_1}\otimes\dots\otimes
g_{n_{i-1}}\otimes h_{n_i}\otimes f_{n_{i+1}}\otimes\dots$$
$$\dots\otimes f_{n_{m+2}})(\eta^{\otimes
(n_1+1)}\otimes\dots\otimes \eta^{\otimes (n_{i-1}+1)}\otimes
\eta^{\otimes (n_i+1)}\otimes \eta^{\otimes
(n_{i+1}+1)}\otimes\dots\otimes \eta^{\otimes (n_{m+2}+1)})=$$
$$=\left(\sum_{m=0}^n\sum_{\,\,\,J_{n-m}}\sum_{i=1}^{m+2}(-1)^{\varrho}\pi_m(g_{n_1}\otimes\dots\otimes
g_{n_{i-1}}\otimes h_{n_i}\otimes f_{n_{i+1}}\otimes\dots\otimes
f_{n_{m+2}})\right)\eta^{\otimes (n+2)},$$ where
$$\alpha=(n_{m+1}+1)n_{m+2}+(n_m+1)(n_{m+1}+n_{m+2})+\dots$$
$$\dots+(n_i+1)(n_{i+1}+\dots+n_{m+2})=\sum_{k=i}^{m+1}(n_k+1)(n_{k+1}+\dots+n_{m+2}),$$
$$\beta=(n_{i-1}+1)((n_i+1)+n_{i+1}+\dots+n_{m+2})\,+$$
$$+\,(n_{i-2}+1)(n_{i-1}+(n_i+1)+n_{i+1}+\dots+n_{m+2})+\dots$$
$$\dots+(n_1+1)(n_2+\dots+n_{i-1}+(n_i+1)+n_{i+1}+\dots+n_{m+2})=$$
$$=\sum_{k=1}^{i-1}(n_k+1)(n_{k+1}+\dots+n_{m+2})+\sum_{s=1}^{i-1}(n_s+1),$$
$$\varrho=m-i+1+\alpha+\beta=m+\sum_{k=1}^{m+1}(n_k+1)(n_{k+1}+\dots+n_{m+2})+\sum_{s=1}^{i-1}n_s.$$
The obtained above equalities follows that the relations $(2.3)$
holds because the map $\eta$ is a isomorphism of degree $(-1)$ of
graded modules.

By using the notion of a homotopy between morphisms of
$\A$-algebras the notion homotopy equivalent $\A$-algebras is
introduced in the usual way.

In conclusion of this paragraph we mention that the homotopy
invariance of the structure of an $\A$-algebra under homotopy
equivalences of the type of an SDR-data of differential modules
was established in \cite{Kad}. \vspace{0.5cm}

\centerline{\bf \S\,3. The tensor functor from the category of
$\A$-algebras}
\centerline{\bf into the category of differential
modules with $\infty$-simplicial faces.} \vspace{0.5cm}

In \cite{Lap4} (see also \cite{Lap2}) it was shown that each
$\A$-algebra $(A,d,\pi_n)$ defines the tensor differential module
with $\infty$-simplicial faces $(T(A),d,\widetilde{\p})$. Let us
consider the construction of this tensor differential module with
$\infty$-simplicial faces.

Given any $\A$-al\-geb\-ra $(A,d,\pi_n)$, consider the tensor
differential module $(T(A),d)$ defined by
$$T(X)=\{T(A)_{n,m}\},\quad T(A)_{n,m}=(A^{\otimes n})_m,\quad
n\geqslant 0,\quad m\geqslant 0,$$ and $d:T(A)_{n,\bu}\to
T(A)_{n,\bu-1}$ is the usual differential in a tensor product.
Consider also the family of maps
$$\widetilde{\p}=\{\p^n_{(i_1,\dots,i_k)}:T(A)_{n,q}\to
T(A)_{n-k,q+k-1}\},$$ $$n\geqslant 0,\quad q\geqslant 0,\quad
1\leqslant k\leqslant n,\quad 0\leqslant i_1<\dots<i_k\leqslant
n,\quad $$ defined by $$\p^n_{(i_1,\dots,i_k)}\!=
\left\{\begin{array}{ll} (-1)^{k(q-1)}1^{\otimes
(j-1)}\otimes\,\pi_{k-1}\otimes 1^{\otimes
(n-k-j)}\,,\,\,\mbox{if}\,\,1\leqslant j\leqslant
n-k&\\\quad\mbox{and}\,\,\, (i_1,\dots,i_k)=(j,j+1,\dots,j+k-1);\\
0,\quad\mbox{otherwise}.&\\
\end{array}\right.\eqno(3.1)$$

In \cite{Lap4} it was proved the following assertion, which
establish the connection between of $\A$-algebras and differential
modules with $\infty$-simplicial faces.

{\bf Theorem 3.1}. Given any $\A$-algebra $(A,d,\pi_n)$, the above
triple $(T(A),d,\widetilde{\p})$ is a differential module with
$\infty$-simplicial faces.~~~$\blacksquare$

Let us show that each morphism of $\A$-algebra $f:(A,d,\pi_n)\to
(A',d,\pi_n)$ defines the morphism
$\widetilde{f}:(T(A),d,\widetilde{\p})\to
(T(A'),d,\widetilde{\p})$ of differential modules with
$\infty$-sim\-pli\-cial faces.

Given any morphism of $\A$-algebras $f:(A,d,\pi_n)\to
(A',d,\pi_n)$, consider the family of module maps
$$\widetilde{f}=\{f^n_{(i_1,\dots,i_k)}:T(A)_{n,q}\to
T(A')_{n-k,q+k}\},$$ $$n\geqslant 0,\quad q\geqslant 0,\quad
0\leqslant k\leqslant n,\quad 0\leqslant i_1<\dots<i_k\leqslant
n,$$ defined by the following rules:

1). If $k=0$, then
$$f^n_{(\,\,)}=\underbrace{f_0\otimes\dots\otimes
f_0}_n.\eqno(3.2)$$

2). If $i_1=0$ or $i_k=n$, then

$$f^n_{(i_1,\dots,i_k)}=0.\eqno(3.3)$$

3). If $i_1\not=0$, $i_k\not=n$,
$(i_1,\dots,i_k)=((j^1_1,\dots,j^1_{n_1}),(j^2_1,\dots,j^2_{n_2}),\dots,(j^s_1,\dots,j^s_{n_s}))$,
$$1\leqslant s\leqslant k,\quad n_1\geqslant 1,\dots,n_s\geqslant
1,\quad n_1+\dots+n_s=k,$$ $$j_{p+1}^m=j_p^m+1,\quad 1\leqslant
p\leqslant n_m-1,\quad 1\leqslant m\leqslant s,\quad j_1^{m+1}>
j^m_{n_m}+1,\quad 1\leqslant m\leqslant s-1,$$ then
$$f^n_{(i_1,\dots,i_k)}=(-1)^{k(q-1)+\gamma}\underbrace{f_0\otimes\dots\otimes
f_0}_{k_1-1}\otimes\,
f_{n_1}\!\otimes\underbrace{f_0\otimes\dots\otimes
f_0}_{k_2-1}\otimes \,f_{n_2}\,\otimes$$
$$\otimes\underbrace{f_0\otimes\dots\otimes
f_0}_{k_3-1}\otimes\,f_{n_3}\!\otimes\dots\otimes\underbrace{f_0\otimes\dots\otimes
f_0}_{k_s-1}\otimes\,f_{n_s}\!\otimes\underbrace{f_0\otimes\dots\otimes
f_0}_{n-(k_1+\dots+k_s)-k},\eqno(3.4)$$ where $$k_1=j_1^1,\quad
k_i=j^{\,i}_1-j^{\,i-1}_{n_{i-1}}-1,\quad 2\leqslant i\leqslant
s,\quad \gamma=\sum\limits_{i=1}^{s-1}n_i(n_{i+1}+\dots+n_s).$$

For example, given the map $f^{15}_{(2,3,6,7,8)}:T(A)_{15,q}\to
T(A')_{10,q+5}$, the formula $(3.4)$ take the following view
$$f^{15}_{((2,3),(6,7,8))}=(-1)^{5(q-1)+2\cdot 3}f_0\otimes
f_2\otimes f_0\otimes f_3\otimes\underbrace{f_0\otimes\dots\otimes
f_0}_6.$$

It is worth mentioning that an arbitrary collection of numbers
$(i_1,\dots,i_k)$, where $1\leqslant i_1<\dots<i_k\leqslant n-1$,
always can be represented in the view, which specified in the
point 3).

The following assertion establish the connection between morphisms
of $\A$-al\-geb\-ras and morphisms of differential modules with
$\infty$-simplicial faces.

{\bf Theorem 3.2}. Given any morphism of $\A$-algebras
$f:(A,d,\pi_n)\to (A',d,\pi_n)$, the family of maps
$\widetilde{f}=\{f^n_{(i_1,\dots,i_k)}:T(A)_{*,\bu}\to
T(A')_{*-k,\bu+k}\}$ defined by $(3.2)$\,--\,$(3.4)$ is the
morphism $\widetilde{f}:(T(A),d,\widetilde{\p})\to
(T(A'),d,\widetilde{\p})$ of differential modules with
$\infty$-sim\-pli\-cial faces.

{\bf Proof}. For the family of maps
$\widetilde{f}=\{f^n_{(i_1,\dots,i_k)}:T(A)_{n,q}\to
T(A)_{n-k,q+k}\}$ defined by $(3.2)$\,--\,$(3.4)$, we need to
check that the relations $(1.2)$ holds. It is clearly that at
$k=0$ we have $d(f^n_{(\,\,)})=0$ because $d(f_0)=0$. It is easy
to see that, for any collection of numbers $(i_1,\dots,i_k)$,
where $i_1=0$, and any permutation $\sigma\in\Sigma_k$, we have
that each partition
$(\widehat{\sigma(i_1)},\dots,\widehat{\sigma(i_m)}\,
|\,\widehat{\sigma(i_{m+1})},\dots,\widehat{\sigma(i_k)})\in
I_\sigma$ satisfies the condition $\widehat{\sigma(i_1)}=0$ or the
condition $\widehat{\sigma(i_{m+1})}=0$. It follows the required
relation $d(f^n_{(0,i_2,\dots,i_k)})=0$. Similarly, for any
collection of numbers $(i_1,\dots,i_k)$, where $i_k=n$, and any
permutation $\sigma\in\Sigma_k$, we have that each partition
$(\widehat{\sigma(i_1)},\dots,\widehat{\sigma(i_m)}\,
|\,\widehat{\sigma(i_{m+1})},\dots,\widehat{\sigma(i_k)})\in
I_\sigma$ satisfies the condition $\widehat{\sigma(i_k)}=n$ or the
condition $\widehat{\sigma(i_m)}=n-(k-m)$. It follows the required
relation $d(f^n_{(i_1,\dots,i_{k-1},n)})=0$. Now, we check that
the maps
$$f^{n+2}_{(1,2,\dots,n+1)}=(-1)^{(n+1)(q-1)}f_{n+1}:(A^{\otimes
(n+2)})_q\to A'_{q+n+1},\quad n\geqslant 0,$$ satisfy the
relations $(1.2)$. First, we rewrite the relations $(2.2)$ in the
form $$d(f_{n+1})=f_0\pi_n
+\sum_{m=0}^{n-1}\sum_{t=1}^{m+2}(-1)^{t(n-m)+n+1}f_{m+1}(\underbrace{1\otimes\dots\otimes
1}_{t-1}\otimes\,\pi_{n-m-1}\otimes\underbrace{1\otimes\dots\otimes
1}_{m+2-t})\,-$$ $$-\,\pi_n(f_0\otimes\dots\otimes f_0)-
\sum_{m=0}^{n-1}\sum_{s=1}^{m+2}\sum_{N_{n-m}}\sum_{T_{m+2}}(-1)^\mu\pi_m(\underbrace{f_0\otimes\dots\otimes
f_0}_{t_1-1}\otimes \,f_{n_1}\otimes$$
$$\otimes\underbrace{f_0\otimes\dots\otimes
f_0}_{t_2-1}\otimes\,f_{n_2}\otimes\dots\otimes\underbrace{f_0\otimes\dots\otimes
f_0}_{t_s-1}\otimes
\,f_{n_s}\!\otimes\underbrace{f_0\otimes\dots\otimes
f_0}_{m+2-(t_1+\dots+t_s)}),\quad n\geqslant -1,\eqno(3.5)$$ where
$$N_{n-m}=\{n_1\geqslant 1,n_2\geqslant 1,\dots,n_s\geqslant
1~|~n_1+n_2+\dots+n_s=n-m\},$$ $$T_{m+2}=\{t_1\geqslant
1,\dots,t_s\geqslant 1~|~t_1+\dots+t_s\leqslant m+2\},$$
$$\mu=\sum_{i=1}^s(t_i-1)(n_i+\dots+n_s)+\sum_{i=1}^{s-1}(n_i+1)(n_{i+1}+\dots+n_s).$$
Given any fixed collections $(n_1,\dots,n_s)\in I_{n-m}$ and
$(t_1,\dots,t_s)\in T_{m+2}$, consider the partition
$(a_1,b_1,a_2,b_2,\dots,a_s,b_s,a_{s+1})$ of the collection of
numbers $(1,2,\dots,n+1)$ into the following $2s+1$ blocks
$$a_1=(1,2,\dots,t_1-1),\quad b_1=(t_1,t_1+1,\dots,t_1+n_1-1),$$
$$a_i=(\sum_{k=1}^{i-1}t_k+n_k,\sum_{k=1}^{i-1}t_k+n_k+1,\dots,\sum_{k=1}^{i-1}t_k+n_k+t_i-1),\quad
2\leqslant i\leqslant s,$$
$$b_i=(\sum_{k=1}^{i-1}t_k+n_k+t_i,\sum_{k=1}^{i-1}t_k+n_k+t_i+1,\dots,\sum_{k=1}^{i-1}t_k+n_k+t_i+n_i-1),\quad
2\leqslant i\leqslant s,$$
$$a_{s+1}=(\sum_{k=1}^st_k+n_k,\sum_{k=1}^st_k+n_k+1,\dots,n+1).$$
Any specified above partition
$(1,2,\dots,n+1)=(a_1,b_1,a_2,b_2,\dots,a_s,b_s,a_{s+1})$ defines
the permutation $\sigma_{n_1,\dots
n_s,t_1,\dots,t_s}\in\Sigma_{n+1}$, which act on the collection of
numbers $(1,2,\dots,n+1)$ by the following rule:
$$\sigma_{n_1,\dots
n_s,t_1,\dots,t_s}(1,2,\dots,n+1)=\sigma_{n_1,\dots
n_s,t_1,\dots,t_s}(a_1,b_1,a_2,b_2,\dots,a_s,b_s,a_{s+1})=$$
$$=(a_1,a_2,\dots,a_s,a_{s+1},b_1,b_2,\dots,b_s).\eqno(3.6)$$ It
is easy to check that the equality of collections
$$(\widehat{\sigma_{n_1,\dots
n_s,t_1,\dots,t_s}(1)},\dots,\widehat{\sigma_{n_1,\dots
n_s,t_1,\dots,t_s}(n+1))}=(1,2,\dots,m+1,b_1,b_2,\dots,b_s)\eqno(3.7)$$
is true. The equalities $(3.1)$\,--\,$(3.4)$ and $(3.7)$ imply
that the relations $(4)$ in the considered case can be rewritten
in the form
$$d(f^{n+2}_{(1,2,\dots,n+1)})=f^1_{(\,\,)}\p^{n+2}_{(1,2,\dots,n+1)}\,+$$
$$+\sum_{m=0}^{n-1} \sum_{t=1}^{m+2}(-1)^{{\rm
sign}(\sigma_{t,n-m})}f^{m+2}_{(1,2,\dots,m+1)}\p^{n+2}_{(t,t+1,\dots,t+n-m-1)}
-\p^{n+2}_{(1,2,\dots,n+1)}f^{n+2}_{(\,\,)}\,-$$
$$-\sum_{m=0}^{n-1}\sum_{s=1}^{m+2}
\sum_{N_{n-m}}\sum_{T_{m+2}}(-1)^{{\rm
sign}(\sigma_{t_1,\dots,t_s,n_1,\dots,n_s})}\p^{m+2}_{(1,2,\dots,m+1)}
f^{n+2}_{(b_1,b_2,\dots,b_s)},\eqno(3.8)$$ where by
$\sigma_{t,n-m}$ we denote the permutation $\sigma_{t_1,n_1}$ such
that $t_1=t$ and $n_1=n-m$. Now, we calculate ${\rm
sign}(\sigma_{t_1,\dots,t_s,n_1,\dots,n_s})$, $1\leqslant
s\leqslant m+2$. Denote by $|a_i|$ the number of elements in the
block $a_i$, $1\leqslant i\leqslant s+1$, and by $|b_j|$ the
number of elements in the block $b_j$, $1\leqslant j\leqslant s$.
Since $\sigma_{t_1,\dots,t_s,n_1,\dots,n_s})$ is the permutation
acting on the collection $(1,2,\dots,n+1)$ by partitioning this
collection into blocks as
$(a_1,b_1,a_2,b_2,\dots,a_s,b_s,a_{s+1})$ and permuting this
blocks by the formula $(3.6)$, the number
$I(\sigma_{t_1,\dots,t_s,n_1,\dots,n_s})$ of inversions of the
permutation $\sigma_{t_1,\dots,t_s,n_1,\dots,n_s})$ is equal
$$I(\sigma_{t_1,\dots,t_s,n_1,\dots,n_s})=$$
$$=|a_2||b_1|+|a_3|(|b_1|+|b_2|)+
\dots+|a_s|(|b_1|+\dots+|b_{s-1}|)+|a_{s+1}|(|b_1|+\dots+|b_s|).$$
Taking into account the congruence
$I(\sigma_{t_1,\dots,t_s,n_1,\dots,n_s})\equiv {\rm
sign}(\sigma_{t_1,\dots,t_s,n_1,\dots,n_s})\,{\rm mod}(2)$, and by
using the equalities $$|a_i|=t_i,\quad |b_i|=n_i,\quad 1\leqslant
i\leqslant s,\quad |a_{s+1}|=n+2-\sum_{k=1}^s(t_k+n_k),\quad
\sum_{i=1}^s n_i=n-m,$$ we obtain the congruence $${\rm
sign}(\sigma_{t_1,\dots,t_s,n_1,\dots,n_s})\equiv
t_2n_1+t_3(n_1+n_2)+\dots+t_s(n_1+\dots+n_{s-1})\,+$$
$$+\,mn+m+(t_1+\dots+t_s)(n-m)\,{\rm mod}(2).$$ In particular, at
$s=1$, $t_1=t$, $n_1=n-m$ we have the following congruence $${\rm
sign}(\sigma_{t,n-m})\equiv mn+m+t(n-m)\,{\rm mod}(2).$$ Now, we
show that the relations $(3.8)$ are equivalence to the relations
$(3.5)$. Indeed, by using the equalities $(3.1)$, $(3.2)$ and
$(3.4)$ we rewrite the relations $(3.8)$ in the form
$$d(f_{n+1})=f_0\pi_n
+\sum_{m=0}^{n-1}\sum_{t=1}^{m+2}(-1)^{\alpha}f_{m+1}(\underbrace{1\otimes\dots\otimes
1}_{t-1}\otimes\,\pi_{n-m-1}\otimes\underbrace{1\otimes\dots\otimes
1}_{m+2-t})\,-$$ $$-\,\pi_n(f_0\otimes\dots\otimes f_0)-
\sum_{m=0}^{n-1}\sum_{s=1}^{m+2}\sum_{N_{n-m}}\sum_{T_{m+2}}(-1)^\beta\pi_m(\underbrace{f_0\otimes\dots\otimes
f_0}_{t_1-1}\otimes \,f_{n_1}\otimes$$
$$\otimes\underbrace{f_0\otimes\dots\otimes
f_0}_{t_2-1}\otimes\,f_{n_2}\otimes\dots\otimes\underbrace{f_0\otimes\dots\otimes
f_0}_{t_s-1}\otimes
\,f_{n_s}\!\otimes\underbrace{f_0\otimes\dots\otimes
f_0}_{m+2-(t_1+\dots+t_s)}),$$ where $$\alpha=(n+1)(q-1)+{\rm
sign}(\sigma_{t,n-m})+(n-m)(q-1)\,+$$ $$+\,(m+1)(q+(n-m-1)-1),$$
$$\beta=(n+1)(q-1)+{\rm
sign}(\sigma_{t_1,\dots,t_s,n_1,\dots,n_s})+(n-m)(q-1)\,+$$
$$+\sum_{i=1}^{s-1}n_i(n_{i+1}+\dots+n_s)+(m+1)(q+(n-m)-1).$$ For
the exponent $\alpha$, we have $$\alpha\equiv
(n+1)(q-1)+mn+m+t(n-m)+(n-m)(q-1)\,+$$
$$+\,(m+1)(q+(n-m-1)-1)\equiv (m+1)(n-m-1)+mn+m+t(n-m)\equiv$$
$$\equiv t(n-m)+n+1\,{\rm mod}(2).$$ For the exponent $\beta$, by
using the equality $n-m=n_1+n_2+\dots+n_s$ we have $$\beta\equiv
(n+1)(q-1)+t_2n_1+t_3(n_1+n_2)+\dots+t_s(n_1+\dots+n_{s-1})\,+$$
$$+\,mn+m+(t_1+\dots+t_s)(n-m)+(n-m)(q-1)\,+$$
$$+\,\sum_{i=1}^{s-1}n_i(n_{i+1}+\dots+n_s)+(m+1)(q+(n-m)-1)\equiv$$
$$\equiv t_2n_1+t_3(n_1+n_2)+\dots+t_s(n_1+\dots+n_{s-1})\,+$$
$$+\,n-m+(t_1+\dots+t_s)(n-m)+\sum_{i=1}^{s-1}n_i(n_{i+1}+\dots+n_s)\equiv$$
$$\equiv t_2n_1+t_3(n_1+n_2)+
\dots+t_s(n_1+\dots+n_{s-1})+(t_1-1)(n_1+\dots+n_s)\,+$$
$$+\,(t_2+\dots+t_s)(n_1+\dots+n_s)+
\sum_{i=1}^{s-1}n_i(n_{i+1}+\dots+n_s)\equiv$$ $$\equiv
\sum_{i=1}^s(t_i-1)(n_i+\dots+n_s)+\sum_{i=1}^{s-1}(n_i+1)(n_{i+1}+\dots+n_s).$$
Thus the relations $(3.8)$ are equivalent to the relations $(3.5)$
and, consequently, the maps
$f^{n+2}_{(1,2,\dots,n+1)}=(-1)^{(n+1)(q-1)}f_{n+1}:(A^{\otimes
(n+2)})_q\to A'_{q+n+1}$, $n\geqslant 0$, satisfy the relations
$(1.2)$. In the same way, as it was done above, it is checked that
all maps $$f^{n+2}_{(j,j+1,\dots,j+k)}=$$
$$=(-1)^{(k+1)(q-1)}\underbrace{f_0\otimes\dots\otimes
f_0}_{j-1}\otimes\,
f_{k+1}\otimes\underbrace{f_0\otimes\dots\otimes f_0}_{n-(j-1)-k}
:(A^{\otimes (n+2)})_q\to (A'^{\otimes(n-k+1)})_{q+k+1},$$
$$n\geqslant 0,\quad 0\leqslant k\leqslant n,\quad1\leqslant
j\leqslant n-k+1,$$ satisfy the relations $(1.2)$. By using this
and also by using the obvious equality
$$d(\underbrace{f_0\otimes\dots\otimes f_0}_{k_1-1}\otimes\,
f_{n_1}\!\otimes\dots\otimes\underbrace{f_0\otimes\dots\otimes
f_0}_{k_s-1}\otimes\,f_{n_s}\!\otimes\underbrace{f_0\otimes\dots\otimes
f_0}_{n-(k_1+\dots+k_s)-k})=$$
$$=\sum_{i=1}^s(-1)^{n_1+\dots+n_{i-1}}\underbrace{f_0\otimes\dots\otimes
f_0}_{k_1-1}\otimes\,
f_{n_1}\!\otimes\dots\otimes\underbrace{f_0\otimes\dots\otimes
f_0}_{k_i-1}\otimes\, d(f_{n_i})\!\otimes\dots$$
$$\dots\otimes\underbrace{f_0\otimes\dots\otimes
f_0}_{k_s-1}\otimes\,f_{n_s}\!\otimes\underbrace{f_0\otimes\dots\otimes
f_0}_{n-(k_1+\dots+k_s)-k},$$ we easily obtain that all maps
$f^n_{(i_1,\dots,i_k)}$ defined by $(3.4)$ satisfy the relations
$(1.2)$.~~~$\blacksquare$

The following assertion establish the connection between the
composition of morphisms of $\A$-al\-geb\-ras and the composition
of morphisms of differential modules with $\infty$-simplicial
faces. This assertion is proved in the similar way as Theorem 2.

{\bf Theorem 3.3}. Given any morphisms of $\A$-algebras
$f:(A,d,\pi_n)\to (A',d,\pi_n)$ and $g:(A',d,\pi_n)\to
(A'',d,\pi_n)$, there is the equality
$$\widetilde{gf}=\widetilde{g}\widetilde{f}:(T(A),d,\widetilde{\p})\to
(T(A''),d,\widetilde{\p})$$ of morphisms of differential modules
with $\infty$-simplicial faces.~~~$\blacksquare$

Theorems $3.1$\,--\,$3.3$ implies that there is the functor of
taking of a tensor module or, more briefly, the tensor functor
$$T:\A(K)\to D\F(K),\quad
T((A,d,\pi_n))=(T(A),d,\widetilde{\p}),$$ $$T(f:(A,d,\pi_n)\to
(A',d,\pi_n))=\widetilde{f}:(T(A),d,\widetilde{\p})\to
(T(A'),d,\widetilde{\p}),$$ where by $\A(K)$ denoted the category
of $\A$-algebras over the ground ring $K$ and by $D\F(K)$ denoted
the category of differential modules with $\infty$-simplicial
faces over the ground ring $K$.

Let us show that the tensor functor $T:\A(K)\to D\F(K)$ sends
homotopies between morphisms into homotopies between morphisms.
More precisely, let us show that each homotopy $h:(A,d,\pi_n)\to
(A',d,\pi_n)$ between morphism of $\A$-algebras
$f,g:(A,d,\pi_n)\to (A',d,\pi_n)$ defines the homotopy
$\widetilde{h}:(T(A),d,\widetilde{\p})\to
(T(A'),d,\widetilde{\p})$ between morphisms
$\widetilde{f},\widetilde{g}:(T(A),d,\widetilde{\p})\to
(T(A'),d,\widetilde{\p})$ of differential modules with
$\infty$-sim\-pli\-cial faces.

Given any homotopy $h:(A,d,\pi_n)\to (A',d,\pi_n)$ between
morphisms of $\A$-al\-geb\-ras $f,g:(A,d,\pi_n)\to (A',d,\pi_n)$,
consider the family of module maps
$$\widetilde{h}=\{h^n_{(i_1,\dots,i_k)}:T(A)_{n,q}\to
T(A')_{n-k,q+k+1}\},$$ $$n\geqslant 0,\quad q\geqslant 0,\quad
0\leqslant k\leqslant n,\quad 0\leqslant i_1<\dots<i_k\leqslant
n,$$ defined by the following rules:

1). If $k=0$, then
$$h^n_{(\,\,)}=\sum\limits_{i=1}^n\underbrace{g_0\otimes\dots\otimes
g_0}_{i-1}\otimes\,h_0\otimes\underbrace{f_0\otimes\dots\otimes
f_0}_{n-i}.\eqno(3.9)$$

2). If $i_1=0$ or $i_k=n$, then

$$h^n_{(i_1,\dots,i_k)}=0.\eqno(3.10)$$

3). If $i_1\not=0$, $i_k\not=n$,
$(i_1,\dots,i_k)=((j^1_1,\dots,j^1_{n_1}),(j^2_1,\dots,j^2_{n_2}),\dots,(j^s_1,\dots,j^s_{n_s}))$,
$$1\leqslant s\leqslant k,\quad n_1\geqslant 1,\dots,n_s\geqslant
1,\quad n_1+\dots+n_s=k,$$ $$j_{p+1}^m=j_p^m+1,\quad 1\leqslant
p\leqslant n_m-1,\quad 1\leqslant m\leqslant s,\quad j_1^{m+1}>
j^m_{n_m}+1,\quad 1\leqslant m\leqslant s-1,$$ then
$$h^n_{(i_1,\dots,i_k)}=(-1)^{k(q-1)+\gamma}\sum_{i=1}^s(-1)^{n_1+\dots+n_{i-1}}\underbrace{g_0\otimes\dots\otimes
g_0}_{k_1-1}\otimes\,g_{n_1}\!\otimes\dots$$
$$\dots\,\otimes\underbrace{g_0\otimes\dots\otimes
g_0}_{k_{i-1}-1}\otimes\,g_{n_{i-1}}\!\otimes\underbrace{g_0\otimes\dots\otimes
g_0}_{k_i-1}\otimes\,h_{n_i}\!\otimes\underbrace{f_0\otimes\dots\otimes
f_0}_{k_{i+1}-1}\otimes\,f_{n_{i+1}}\!\otimes\dots$$
$$\dots\otimes\underbrace{f_0\otimes\dots\otimes
f_0}_{k_s-1}\otimes\,f_{n_s}\!\otimes\underbrace{f_0\otimes\dots\otimes
f_0}_{k_{s+1}-1}\,+$$
$$+\,(-1)^{k(q-1)+\gamma}\sum_{i=1}^{s+1}(-1)^{n_1+\dots+n_{i-1}}\underbrace{g_0\otimes\dots\otimes
g_0}_{k_1-1}\otimes\,g_{n_1}\!\otimes\dots\otimes\underbrace{g_0\otimes\dots\otimes
g_0}_{k_{i-1}-1}\otimes\,g_{n_{i-1}}\otimes$$
$$\otimes\left\{\sum_{j=1}^{k_i-1}\underbrace{g_0\otimes\dots\otimes
g_0}_{j-1}\otimes\,h_0\otimes\underbrace{f_0\otimes\dots\otimes
f_0}_{k_i-1-j}\right\}\otimes\,f_{n_i}\!\otimes\underbrace{f_0\otimes\dots\otimes
f_0}_{k_{i+1}-1}\otimes\,f_{n_{i+1}}\otimes\dots$$
$$\dots\otimes\underbrace{f_0\otimes\dots\otimes
f_0}_{k_s-1}\otimes\,f_{n_s}\!\otimes\underbrace{f_0\otimes\dots\otimes
f_0}_{k_{s+1}-1}\,,\eqno(3.11)$$ where $$k_1=j_1^1,\quad
k_i=j^{\,i}_1-j^{\,i-1}_{n_{i-1}}-1,\quad 2\leqslant i\leqslant
s,\quad k_{s+1}=n-(k_1+\dots+k_s)-k+1,$$
$$\gamma=\sum\limits_{i=1}^{s-1}n_i(n_{i+1}+\dots+n_s).$$

For example, given the map $h^{15}_{(2,3,6,7,8)}:T(A)_{15,q}\to
T(A')_{10,q+6}$, the formula $(3.11)$ take the following view
$$h^{15}_{(2,3,6,7,8)}=h^{15}_{((2,3),(6,7,8))}=(-1)^{5(q-1)+2\cdot
3}g_0\otimes h_2\otimes f_0\otimes
f_3\otimes\underbrace{f_0\otimes\dots\otimes f_0}_6+$$
$$+(-1)^{5(q-1)+2\cdot 3+2}g_0\otimes g_2\otimes g_0\otimes
h_3\otimes\underbrace{f_0\otimes\dots\otimes f_0}_6+$$
$$+(-1)^{5(q-1)+2\cdot 3}h_0\otimes f_2\otimes f_0\otimes
f_3\otimes\underbrace{f_0\otimes\dots\otimes f_0}_6+$$
$$+(-1)^{5(q-1)+2\cdot 3+2}g_0\otimes g_2\otimes h_0\otimes
f_3\otimes\underbrace{f_0\otimes\dots\otimes f_0}_6+$$
$$+(-1)^{5(q-1)+2\cdot 3+2+3}g_0\otimes g_2\otimes g_0\otimes
g_3\otimes\left\{\sum_{j=1}^6\underbrace{g_0\otimes\dots\otimes
g_0}_{j-1}\otimes\,h_0\otimes\underbrace{f_0\otimes\dots\otimes
f_0}_{6-j}\right\}.$$

As above we must note that an arbitrary collection of numbers
$(i_1,\dots,i_k)$, where $1\leqslant i_1<\dots<i_k\leqslant n-1$,
always can be represented in the view, which specified in the
point 3).

{\bf Theorem 3.4}. Suppose given any homotopy $h:(A,d,\pi_n)\to
(A',d,\pi_n)$ between morphisms of $\A$-algebras
$f,g:(A,d,\pi_n)\to (A',d,\pi_n)$. Then the family of maps
$\widetilde{h}=\{h^n_{(i_1,\dots,i_k)}:T(A)_{*,\bu}\to
T(A')_{*-k,\bu+k+1}\}$ defined by $(3.9)$\,--\,$(3.11)$ is the
homotopy $\widetilde{h}:(T(A),d,\widetilde{\p})\to
(T(A'),d,\widetilde{\p})$ between morphisms
$\widetilde{f},\widetilde{g}:(T(A),d,\widetilde{\p})\to
(T(A'),d,\widetilde{\p})$ of differential modules with
$\infty$-simplicial faces.

{\bf Proof}. For the family of maps
$\widetilde{h}=\{h^n_{(i_1,\dots,i_k)}:T(A)_{n,q}\to
T(A)_{n-k,q+k+1}\}$ defined by $(3.9)$\,--\,$(3.11)$, we need to
check that the relations $(1.3)$ holds. Clearly, that at $k=0$ we
have $d(h^n_{(\,\,)})=f^n_{(\,\,)}-g^n_{(\,\,)}$ because
$d(h_0)=f_0-g_0$. It is easy to see that, for any collection of
numbers $(i_1,\dots,i_k)$, where $i_1=0$, and any permutation
$\sigma\in\Sigma_k$, we have that each partition
$(\widehat{\sigma(i_1)},\dots,\widehat{\sigma(i_m)}\,
|\,\widehat{\sigma(i_{m+1})},\dots,\widehat{\sigma(i_k)})\in
I_\sigma$ satisfies the condition $\widehat{\sigma(i_1)}=0$ or the
condition $\widehat{\sigma(i_{m+1})}=0$. It follows the required
relation $d(h^n_{(0,i_2,\dots,i_k)})=0$. Similarly, for any
collection of numbers $(i_1,\dots,i_k)$, where $i_k=n$, and any
permutation $\sigma\in\Sigma_k$, we have that each partition
$(\widehat{\sigma(i_1)},\dots,\widehat{\sigma(i_m)}\,
|\,\widehat{\sigma(i_{m+1})},\dots,\widehat{\sigma(i_k)})\in
I_\sigma$ satisfies the condition $\widehat{\sigma(i_k)}=n$ or the
condition $\widehat{\sigma(i_m)}=n-(k-m)$. It follows the required
relation $d(h^n_{(i_1,\dots,i_{k-1},n)})=0$. Now, we check that
the maps
$$h^{n+2}_{(1,2,\dots,n+1)}=(-1)^{(n+1)(q-1)}h_{n+1}:(A^{\otimes
(n+2)})_q\to A'_{q+n+2},\quad n\geqslant 0,$$ satisfy the
relations $(1.3)$. First, we rewrite the relations $(2.3)$ in the
form $$d(h_{n+1})=f_{n+1}-g_{n+1}-h_0\pi_n\,+$$
$$+\sum_{m=0}^{n-1}\sum_{t=1}^{m+2}(-1)^{t(n-m)+n}h_{m+1}(\underbrace{1\otimes\dots\otimes
1}_{t-1}\otimes\,\pi_{n-m-1}\otimes\underbrace{1\otimes\dots\otimes
1}_{m-t+2}\,)\,+$$
$$+\sum_{i=1}^{n+2}(-1)^n\pi_n(\underbrace{g_0\otimes\dots\otimes
g_0}_{i-1}\otimes\,h_0\otimes \underbrace{f_0\otimes\dots\otimes
f_0}_{n+2-i}\,)\,+$$
$$+\sum_{m=0}^{n-1}\sum_{s=1}^{m+2}\sum_{N_{n-m}}\sum_{T_{m+2}}\sum_{i=1}^s(-1)^{\vartheta+n_1+\dots+n_{i-1}}
\pi_m(\underbrace{g_0\otimes\dots\otimes g_0}_{t_1-1}\otimes
\,g_{n_1}\!\otimes\dots$$
$$\dots\otimes\underbrace{g_0\otimes\dots\otimes
g_0}_{t_{i-1}-1}\otimes\,g_{n_{i-1}}\!\otimes\underbrace{g_0\otimes\dots\otimes
g_0}_{t_i-1}\otimes\,h_{n_i}\!\otimes\underbrace{f_0\otimes\dots\otimes
f_0}_{t_{i+1}-1}\otimes\,f_{n_{i+1}}\!\otimes\dots$$
$$\dots\otimes\underbrace{f_0\otimes\dots\otimes
f_0}_{t_s-1}\otimes
\,f_{n_s}\!\otimes\underbrace{f_0\otimes\dots\otimes
f_0}_{t_{s+1}-1})\,+$$
$$+\sum_{m=0}^{n-1}\sum_{s=1}^{m+2}\sum_{N_{n-m}}\sum_{T_{m+2}}\sum_{i=1}^{s+1}(-1)^{\vartheta+n_1+\dots+n_{i-1}}
\pi_m(\underbrace{g_0\otimes\dots\otimes g_0}_{t_1-1}\otimes
\,g_{n_1}\!\otimes\dots$$
$$\dots\otimes\underbrace{g_0\otimes\dots\otimes
g_0}_{t_{i-1}-1}\otimes\,g_{n_{i-1}}\otimes\left\{\sum_{j=1}^{t_i-1}\underbrace{g_0\otimes\dots\otimes
g_0}_{j-1}\otimes\,h_0\otimes\underbrace{f_0\otimes\dots\otimes
f_0}_{t_i-1-j}\right\}\otimes f_{n_i}\,\otimes$$
$$\otimes\underbrace{f_0\otimes\dots\otimes
f_0}_{t_{i+1}-1}\otimes\,f_{n_{i+1}}\!\otimes\dots\otimes\underbrace{f_0\otimes\dots\otimes
f_0}_{t_s-1}\otimes
\,f_{n_s}\!\otimes\underbrace{f_0\otimes\dots\otimes
f_0}_{t_{s+1}-1}),\quad n\geqslant -1,\eqno(3.12)$$ where
$$N_{n-m}=\{n_1\geqslant 1,n_2\geqslant 1,\dots,n_s\geqslant
1~|~n_1+n_2+\dots+n_s=n-m\},$$ $$T_{m+2}=\{t_1\geqslant
1,\dots,t_s\geqslant 1~|~t_1+\dots+t_s\leqslant m+2\},\quad
t_{s+1}=m+3-(t_1+\dots+t_s),$$
$$\vartheta=m+\sum_{k=1}^s(t_k-1)(n_k+\dots+n_s)+
\sum_{k=1}^{s-1}(n_k+1)(n_{k+1}+\dots+n_s).$$ The equalities
$(3.1)$\,--\,$(3.4)$, $(3.7)$ and $(3.9)$\,--\,$(3.11)$ imply that
the relations $(1.3)$ in the considered case can be rewritten in
the form
$$d(h^{n+2}_{(1,2,\dots,n+1)})=f^{n+2}_{(1,2,\dots,n+1)}-g^{n+2}_{(1,2,\dots,n+1)}-
h^1_{(\,\,)}\p^{n+2}_{(1,2,\dots,n+1)}\,-$$ $$-\sum_{m=0}^{n-1}
\sum_{t=1}^{m+2}(-1)^{{\rm
sign}(\sigma_{t,n-m})}h^{m+2}_{(1,2,\dots,m+1)}\p^{n+2}_{(t,t+1,\dots,t+n-m-1)}-
\p^{n+2}_{(1,2,\dots,n+1)}h^{n+2}_{(\,\,)}\,-$$
$$-\sum_{m=0}^{n-1}\sum_{s=1}^{m+2}\sum_{N_{n-m}}\sum_{T_{m+2}}(-1)^{{\rm
sign}(\sigma_{t_1,\dots,t_s,n_1,\dots,n_s})}\p^{m+2}_{(1,2,\dots,m+1)}h^{n+2}_{(b_1,b_2,\dots,b_s)},\eqno(3.13)$$
where $\sigma_{t_1,\dots,t_s,n_1,\dots,n_s},
\sigma_{t,n-m}\in\Sigma_{n+1}$ and $b_1,\dots,b_s$ are the same as
in $(3.8)$.

Let us show that the relations $(3.13)$ are equivalent to the
relations $(3.12)$. By using the relations $(3.1)$\,--\,$(3.4)$
and $(3.9)$\,--\,$(3.11)$ rewritten the relation $(3.13)$ in the
form $$d(h_{n+1})=f_{n+1}-g_{n+1}-h_0\pi_n\,+$$
$$+\sum_{m=0}^{n-1}\sum_{t=1}^{m+2}(-1)^{\alpha}h_{m+1}(\underbrace{1\otimes\dots\otimes
1}_{t-1}\otimes\,\pi_{n-m-1}\otimes\underbrace{1\otimes\dots\otimes
1}_{m-t+2}\,)\,+$$
$$+\sum_{i=1}^{n+2}(-1)^n\pi_n(\underbrace{g_0\otimes\dots\otimes
g_0}_{i-1}\otimes\,h_0\otimes \underbrace{f_0\otimes\dots\otimes
f_0}_{n+2-i}\,)\,+$$
$$+\sum_{m=0}^{n-1}\sum_{s=1}^{m+2}\sum_{N_{n-m}}\sum_{T_{m+2}}\sum_{i=1}^s(-1)^{\beta+n_1+\dots+n_{i-1}}
\pi_m(\underbrace{g_0\otimes\dots\otimes g_0}_{t_1-1}\otimes
\,g_{n_1}\!\otimes\dots$$
$$\dots\otimes\underbrace{g_0\otimes\dots\otimes
g_0}_{t_{i-1}-1}\otimes\,g_{n_{i-1}}\!\otimes\underbrace{g_0\otimes\dots\otimes
g_0}_{t_i-1}\otimes\,h_{n_i}\!\otimes\underbrace{f_0\otimes\dots\otimes
f_0}_{t_{i+1}-1}\otimes\,f_{n_{i+1}}\!\otimes\dots$$
$$\dots\otimes\underbrace{f_0\otimes\dots\otimes
f_0}_{t_s-1}\otimes
\,f_{n_s}\!\otimes\underbrace{f_0\otimes\dots\otimes
f_0}_{t_{s+1}-1})\,+$$
$$+\sum_{m=0}^{n-1}\sum_{s=1}^{m+2}\sum_{N_{n-m}}\sum_{T_{m+2}}\sum_{i=1}^{s+1}(-1)^{\beta+n_1+\dots+n_{i-1}}
\pi_m(\underbrace{g_0\otimes\dots\otimes g_0}_{t_1-1}\otimes
\,g_{n_1}\!\otimes\dots$$
$$\dots\otimes\underbrace{g_0\otimes\dots\otimes
g_0}_{t_{i-1}-1}\otimes\,g_{n_{i-1}}\otimes\left\{\sum_{j=1}^{t_i-1}\underbrace{g_0\otimes\dots\otimes
g_0}_{j-1}\otimes\,h_0\otimes\underbrace{f_0\otimes\dots\otimes
f_0}_{t_i-1-j}\right\}\otimes$$ $$\otimes\,
f_{n_i}\!\otimes\underbrace{f_0\otimes\dots\otimes
f_0}_{t_{i+1}-1}\otimes\,f_{n_{i+1}}\!\otimes\dots\otimes\underbrace{f_0\otimes\dots\otimes
f_0}_{t_s-1}\otimes
\,f_{n_s}\!\otimes\underbrace{f_0\otimes\dots\otimes
f_0}_{t_{s+1}-1}),$$ where $$\alpha=(n+1)(q-1)+{\rm
sign}(\sigma_{t,n-m})+(n-m)(q-1)\,+$$ $$+\,(m+1)(q+(n-m-1)-1)+1,$$
$$\beta=(n+1)(q-1)+{\rm
sign}(\sigma_{t_1,\dots,t_s,n_1,\dots,n_s})+(n-m)(q-1)+$$
$$+\sum_{i=1}^{s-1}n_i(n_{i+1}+\dots+n_s)+(m+1)(q+(n-m+1)-1)+1.$$
For the exponent $\alpha$, by using ${\rm
sign}(\sigma_{t,n-m})\equiv mn+m+t(n-m)\,{\rm mod}(2)$ we have
$$\alpha\equiv (n+1)(q-1)+mn+m+t(n-m)+(n-m)(q-1)\,+$$
$$+\,(m+1)(q+(n-m-1)-1)+1\equiv (m+1)(n-m-1)+mn+m+t(n-m)+1\equiv$$
$$\equiv t(n-m)+n\,{\rm mod}(2).$$ For the exponent $\beta$, by
using $${\rm sign}(\sigma_{t_1,\dots,t_s,n_1,\dots,n_s})\equiv
t_2n_1+t_3(n_1+n_2)+\dots+t_s(n_1+\dots+n_{s-1})\,+$$
$$+\,mn+m+(t_1+\dots+t_s)(n-m)\,{\rm mod}(2)$$ and
$n-m=n_1+n_2+\dots+n_s$, we have $$\beta\equiv
(n+1)(q-1)+t_2n_1+t_3(n_1+n_2)+\dots+t_s(n_1+\dots+n_{s-1})\,+$$
$$+\,mn+m+(t_1+\dots+t_s)(n-m)+(n-m)(q-1)+\sum_{i=1}^{s-1}n_i(n_{i+1}+\dots+n_s)\,+$$
$$+\,(m+1)(q+n-m)+1\equiv
t_2n_1+t_3(n_1+n_2)+\dots+t_s(n_1+\dots+n_{s-1})-n\,+$$
$$+\,(t_1+\dots+t_s)(n-m)+\sum_{i=1}^{s-1}n_i(n_{i+1}+\dots+n_s)\equiv$$
$$\equiv
t_2n_1+t_3(n_1+n_2)+\dots+t_s(n_1+\dots+n_{s-1})+m+(t_1-1)(n_1+\dots+n_s)\,+$$
$$+\,(t_2+\dots+t_s)(n_1+\dots+n_s)+\sum_{i=1}^{s-1}n_i(n_{i+1}+\dots+n_s)\equiv$$
$$\equiv m+\sum_{k=1}^s(t_k-1)(n_k+\dots+n_s)+
\sum_{k=1}^{s-1}(n_k+1)(n_{k+1}+\dots+n_s)\,{\rm mod}(2).$$ Thus
the relations $(3.13)$ are equivalent to the relations $(3.12)$
and, consequently, the maps
$h^{n+2}_{(1,2,\dots,n+1)}=(-1)^{(n+1)(q-1)}h_{n+1}:(A^{\otimes
(n+2)})_q\to A'_{q+n+2}$, $n\geqslant 0$, satisfy the relations
$(1.3)$. In the same way, as it was done above, it is checked that
all maps $$h^{n+2}_{(j,j+1,\dots,j+k)}=$$
$$=(-1)^{(k+1)(q-1)}\underbrace{g_0\otimes\dots\otimes
g_0}_{j-1}\otimes\,
h_{k+1}\otimes\underbrace{f_0\otimes\dots\otimes f_0}_{n-(j-1)-k}
:(A^{\otimes (n+2)})_q\to (A'^{\otimes(n-k+1)})_{q+k+2},$$
$$n\geqslant 0,\quad 0\leqslant k\leqslant n,\quad1\leqslant
j\leqslant n-k+1,$$ satisfy the relations $(1.3)$. By using this
and also by using the usual formula of a differential in tensor
products we easily obtain that all maps $h^n_{(i_1,\dots,i_k)}$
defined by $(3.11)$ satisfy the relations
$(1.3)$.~~~$\blacksquare$

By using Theorem 3.4 we easily obtain the following assertion.

{\bf Corollary 3.1}. The tensor functor $T:\A(K)\to D\F(K)$ sends
homotopy equivalent $\A$-algebras into homotopy equivalent
differential modules with $\infty$-sim\-pli\-cial
faces.~~~$\blacksquare$

\vspace{1cm}

Serov Str., Saransk, Russia,

e-mail: slapin@mail.ru


\begin{thebibliography}{99}
\bibitem{Lap2} S. V. Lapin, Homotopy simplicial faces and the homology of realizations of simplicial
topological spaces (in Russian), Mat. Zametki 94 (5), 661–681
(2013); translation in Math. Notes 94 (5–6), 619–635 (2013).
\bibitem{Lap1} S. V. Lapin, Homotopy properties of differential Lie modules over curved
coalgebras and Koszul duality (in Russian), Mat. Zametki 94 (3),
354–372 (2013); translation in Math. Notes 94 (3–4), 335–350
(2013).
\bibitem{Lap3} S. V. Lapin, Differential Lie modules over curved colored coalgebras
and $\infty$-sim\-pli\-cial modules (in Russian), Mat. Zametki 96
(5), 709–731 (2014); translation in Math. Notes 96 (5–6), 698–715
(2014).
\bibitem{Lap4} S. V. Lapin, Chain realization of differential modules with $\infty$-simplicial
faces and the $B$-construction over $\A$-algebras (in Russian),
Mat. Zametki 98 (1), 101–124 (2015); translation in Math. Notes 98
(1–2), 111–129 (2015).
\bibitem{Lap5} S. V. Lapin, Homotopy properties of differential modules
with simplicial $F_\infty$-faces and $\D$-differential modules,
Georgian Math. J. 2015; 22(4):543-562.
\bibitem{Lap6} S. V. Lapin, Homotopy properties of $\infty$-simplicial coalgebras and homotopy unital supplemented
$\A$-algebras (in Russian), Mat. Zametki 99 (1), 55–77 (2016);
translation in Math. Notes 99 (1–2), 63–81 (2016).
\bibitem{Lap7} S. V. Lapin, Differential modules with $\infty$-simplicial faces and
$A_\infty$-algebras, arXiv:1809.01853v1 [math.AT] 6 Sep 2018, p.
1-26.
\bibitem{Lap9} S. V. Lapin, Differential perturbations and $\D$-differential modules (in Russian),
Mat. Sb. 192:11 (2001), 55–76; translation in Sb. Math. 192:11
(2001), 1639–1659.
\bibitem{Lap10} S. V. Lapin, $\D$-differential $\A$-algebras and spectral sequences
(in Russian), Mat. Sb. 193:1 (2002), 119–142; translation in Sb.
Math. 193:1 (2002), 119–142.
\bibitem{Lap11} S. V. Lapin, $(DA)_\infty$-modules over $(DA)_\infty$-algebras and spectral sequences
(in Russian), Izv. RAN Ser. Math. , 66:3 (2002), 103–130;
translation in Izv. Math., 66:3 (2002), 543–568.
\bibitem{Lap12} S. V. Lapin,  $\D$-differentials and
$\A$-structures in spectral sequences, J. Math. Sci., 123:4
(2004), 4221–4254.
\bibitem{Lap13} S. V. Lapin, $\D$-differential $\E$-algebras and multiplicative spectral sequences
(in Russian), Mat. Sb. 196:11 (2005), 75–108; translation in Sb.
Math. 196:11 (2005), 1627–1658.
\bibitem{Lap14} S. V. Lapin, $\D$-differential $\E$-algebras and spectral sequences of fibrations
(in Russian), Mat. Sb. 198:10 (2007), 3-30; translation in Sb.
Math. 198:10 (2007), 1379–1406.
\bibitem{Lap15} S. V. Lapin, $\D$-differential $\E$-algebras and
Steenrod operations in spectral sequences, J. Math. Sci., 152:3
(2008), 372–403.
\bibitem{Lap16} S. V. Lapin, $\D$-differential $\E$-algebras and
spectral sequences of $\D$-dif\-fe\-ren\-tial modules, J. Math.
Sci., 159:6 (2009), 819–832.
\bibitem{Lap17} S. V. Lapin,  Multiplicative $\A$-structure in terms of spectral sequences of fibrations
(in Russian), Fundam. Prikl. Mat., 14:6 (2008), 141–175;
translation in J. Math. Sci., 164:1 (2010), 95–118.
\bibitem{Lap18} S. V. Lapin, Cyclic modules with $\infty$-simplicial faces and cyclic homology of $\A$-al\-geb\-ras
(in Russian), Mat. Zametki 102 (6), 859-880 (2017); translation in
Math. Notes 102, (5–6), 806–823 (2017).
\bibitem{Lap19} S. V. Lapin, Cyclic homology of cyclic $\infty$-simplicial
modules, Georgian Math. J. 2018; 25(4):571–591.
\bibitem{Lap20} S. V. Lapin, Dihedral and reflexive modules with $\infty$-simplicial
faces and dihedral and reflexive homology of involutive
$\A$-algebras over unital commutative rings, ArXiv:1809.07510v1
[math.AT] 20 Sep 2018, p. 1-25.
\bibitem{Kad}  T. V. Kadeishvili, On the theory of homology of fiber spaces (in Russian),
Uspekhi Mat. Nauk 35 (1980), no. 3, 183–188; translation in
Russian Math. Surveys 35 (1980), no. 3, 231–238.
\bibitem{S} J. D. Stasheff, Homotopy associativity of $H$-spaces I,II,
Trans. Amer. Math. Soc., 108:2 (1963), 275-312.
\bibitem{Smir} V. A. Smirnov, Simplicial and Operad Methods in Algebraic
Topology, Transl. Math. Monogr. 198, Amer. Math. Soc., Providence,
2001.

\end{thebibliography}
\end{document}